\documentclass[12pt,a4paper]{article}

\usepackage[english,francais]{babel} 
\RequirePackage[latin1]{inputenc}
\RequirePackage[T1]{fontenc}
\RequirePackage{graphicx,hyperref,ifpdf,color}

\RequirePackage{amsmath,amssymb,amsfonts,bbm}

\newcommand{\be}{{\bf e}}

\newcommand{\bs}{{\bf s}}

\newcommand{\bP}{{\bf P}}
\newcommand{\bE}{{\bf E}}
\newcommand{\bN}{{\bf N}}

\newcommand{\R}{\mathbb{R}}
\newcommand{\N}{\mathbb{N}}

\newcommand{\D}{\mathbb{D}}

\newcommand{\TT} { {\cal T }}

\def\build#1_#2^#3{\mathrel{
\mathop{\kern 0pt#1}\limits_{#2}^{#3}}}

\newcommand{\HT}{{\rm ht}}

\def\cq{$\hfill \square$}
\def\un{\underline}
\def\d{{\rm d}}

\def\eps{\varepsilon}

\def\ba{\begin{eqnarray*}}
\def\ea{\end{eqnarray*}}
\def\ov{\overline}
\def\wh{\widehat}
\def\wt{\widetilde}
\def\xg{\underleftarrow{X}}
\def\xd{\underrightarrow{X}}

\def\rem{\noindent{\bf Remark. }}

\def\proof{\noindent{\bf Proof. }}

\newcommand{\ind}{\mathbbm{1}}

\newtheorem{thm}{Theorem}
\newtheorem{lmm}{Lemma}
\newtheorem{prp}{Proposition}
\newtheorem{defn}{Definition}

\usepackage{hyperref} 
\usepackage{verbatim,multicol}

\def\miermont{{\href{http://www.dma.ens.fr/~miermont/}{Gr\'egory Miermont}}}
\def\dma{{\href{http://www.dma.ens.fr/}{DMA}}}
\def\lpma{{\href{http://www.proba.jussieu.fr/}{LPMA}}}

\begin{document}
\title{Self-similar fragmentations derived 
from the stable tree II: 
splitting at nodes}
\author{\miermont
\\ \dma, \'Ecole Normale Sup\'erieure, \\
and \lpma, Universit\'e Paris VI.\\
45, rue d'Ulm, \\
75230 Paris Cedex 05}
\date{}

\maketitle

\begin{abstract}
We study a natural fragmentation process of the so-called stable tree 
introduced by Duquesne and Le Gall, which consists in removing the nodes 
of the tree according to a certain procedure that makes the fragmentation 
self-similar  with positive  index. Explicit formulas for the semigroup are
given, and we provide asymptotic results. 
We also give an alternative construction of this fragmentation, using paths of
L\'evy processes, hence echoing the two alternative
constructions of the standard additive coalescent by
fragmenting the Brownian continuum random  tree or using
Brownian paths, respectively due to Aldous-Pitman and Bertoin.
\end{abstract}


\bigskip

\noindent{\bf A.M.S. Classification. }60J25, 60G52.

\newpage

\section{Introduction}\label{intro}

The goal of this paper is to investigate a Markovian fragmentation 
of the so-called 
{\em stable tree}. It is a model of {\em continuum random tree} (CRT) 
depending on a parameter $\alpha\in(1,2]$ that has 
been introduced recently by Duquesne and Le Gall \cite{duqleg02}, and which
basically corresponds to a possible scaling limit as $n\to\infty$ of a size $n$
Galton-Watson tree with given progeny distribution. The stable tree is denoted by $\TT$. It
is a random metric space with distance $d$, whose elements $v$ are called 
{\em vertices}. One of these vertices is distinguished and called the 
{\em root}. This space is a tree in that for $v,w\in\TT$, 
there is a unique non-self-crossing path $[[v,w]]$ from $v$ to $w$ in 
$\TT$, whose length equals $d(v,w)$. 
For every $v\in \TT$, call {\em height} of $v$ in $\TT$ and denote by 
$\HT(v)$ the distance of $v$ to the root. 
The {\em leaves} ${\cal L}(\TT)$ of $\TT$ are those vertices that do not belong to 
the interior of any 
path leading from one vertex to another, and the {\em skeleton} of the 
tree is the set $\TT\setminus{\cal L}(\TT)$ of non-leaf vertices. The {\em branchpoints} 
are the vertices $b$ so that there exist $v\neq b,w\neq b$ such that 
$[[{\rm root},v]]\cap[[{\rm root},w]]=[[{\rm root}, b]]$.
With each realization of $\TT$
is associated the uniform probability measure $\mu$, 
called the {\em mass measure}, that is supported by ${\cal L}(\TT)$.  
Details are given in Section \ref{levytree}. 

When $\alpha=2$, the stable tree is, up to a scale factor, 
the Brownian CRT of Aldous \cite{aldouscrt93}.  It has been shown by Aldous
and Pitman \cite{jpda98sac} that a certain device for logging this tree gives rise
to a {\em fragmentation} process which is the time-reversed process of the
so-called {\em standard additive coalescent}.  The idea is as follows. The
Brownian CRT $\TT$ is described by a $\sigma$-finite {\em length measure} $\ell$
carried  by the  skeleton  (non-leaf vertices),  and  a (uniform)  probability
measure $\mu$ on its leaves, called the mass measure. For $t\geq 0$, 
consider a Poisson random measure on $\TT$ with intensity $t\ell$, in a
consistent way as  $t$ varies. When the marked vertices  of the tree are
removed, the tree is decomposed into a random forest, whose ranked
$\mu$-masses form an element $F_{\rm AP}(t)$ of the space 
$$S:=\left\{\bs=(s_1,s_2,\ldots): s_1\geq s_2\geq\ldots\geq 0, 
\sum_{i=1}^{\infty}s_i\leq 1\right\}.$$ 
It is  actually checked that the sum  of components of $F_{\rm  AP}(t)$ is $1$
a.s. Then Bertoin \cite{bertsfrag02} noticed (it was implicit in 
\cite{jpda98sac}) that the process  $(F_{\rm AP}(t),t\geq 0)$ is an $S$-valued
self-similar fragmentation with index $1/2$, in the following sense. 

\begin{defn}
An $S$-valued self-similar fragmentation with index $\beta\in\R$ is an
$S$-valued  Markov  process  starting  a.s.\  from  $(1,0,\ldots)$,  which  is
continuous in probability and satisfies the following 
{\em fragmentation property}:
\begin{quote}
Given  $F(t)=\bs=(s_1,s_2,\ldots)$,  the  law  of  $F(t+t')$ is  that  of  the
decreasing  rearrangement of  the  sequences $s_iF^{(i)}(s_i^{\beta}t'),i\geq
1$, where the $F^{(i)}$'s are independent copies of $F$. 
\end{quote}
\end{defn}

Such fragmentations have been introduced and extensively studied by Bertoin in
\cite{berthfrag01,bertsfrag02}. By \cite{berest02}, the laws of the 
self-similar fragmentations are characterized by a $3$-tuple $(\beta,c,\nu)$, 
where $\beta$ is the self-similarity index, $c\geq 0$ is an 
erosion coefficient and, more importantly, $\nu$ is a $\sigma$-finite 
{\em dislocation measure} on $S$ that integrates the map $\bs\mapsto 1-s_1$. 
This measure $\nu$ describes the ``jumps'' of the fragmentation process, 
i.e.\ the way sudden dislocations occur. Roughly speaking, 
$x^{\beta}\nu(\d\bs)$ is the  instantaneous rate at which an  object with size
$x$ fragments  to form objects with  sizes $x\bs$ (see  also Lemma \ref{msprp}
below). In \cite{bertsfrag02},
Bertoin showed that the erosion of $F_{\rm AP}$ is $0$, and that the
dislocation measure $\nu_{\rm AP}$ is characterized by the two formulas
$$\nu_{\rm AP}(s_1\in \d x)=\frac{\d x}{\sqrt{2\pi x^3(1-x)^3}}\, ,\qquad x\in[1/2,1),$$
and $\nu_{\rm AP}\{\bs: s_1+s_2<1\}=0$ (such fragmentations are called {\em binary}). 

The main motivation of the present paper is to seek for a possible
generalization of  the fragmentation  $F_{\rm AP}$, when  the Brownian  CRT is
replaced by the general $\alpha\in(1,2)$-stable tree. The game is made
interesting in that there are
important structural differences between the Brownian tree and the other stable
trees, which imply that the Aldous-Pitman fragmentation device explained above
(homogeneous fragmentation on the skeleton) gives rise to a
binary fragmentation process which is {\em not} self-similar. 
It seems that the fragmentations hence obtained
are related to the ones studied in \cite{mier01} in relation with the additive
coalescent, but this will be studied elsewhere. The defect in the
self-similarity property comes from the fact that, contrary to the Brownian 
tree which is {\em binary} (its branchpoints have degree $3$), 
the branchpoints of the stable tree are {\em hubs} with infinite 
degree and with different ``magnitudes''. These are not
affected by the Aldous-Pitman fragmentation device, which a.s.\ never cuts at 
branchpoints. Therefore, as time passes,
this device creates small trees with unusually ``large'' hubs, which cannot be
rescaled copies of the initial stable tree. Rather, to obtain self-similarity,
it is needed to directly remove the hubs themselves with a certain strategy. 

Call ${\cal H}(\TT)$ the set of branchpoints of $\TT$, which will also be
referred to as the set of hubs of $\TT$ when dealing with the stable
($\alpha\in(1,2)$) tree. 
To evaluate the magnitude of 
$b\in{\cal H}(\TT)$, consider the {\em fringe subtree} $\TT_b$ rooted at $b$, 
i.e.\ the subset $\{v\in \TT: b\in [[{\rm root},v]]\}$. 
Then one can define the {\em local time}, or {\em width} of the hub $b$ 
as the limit 
\begin{equation}\label{widthhub}
L(b)=\lim_{\eps\downarrow0}\frac{1}{\eps}\mu\{v\in \TT_b:d(v,b)<\eps\}
\end{equation}
which exists a.s.\ and is positive: see Proposition \ref{onetoone} below. 

Now given a realization of $\TT$ and for every $b\in{\cal H}(\TT)$,
take a standard exponential random 
variable $e_b$, so that the variables $e_b$ are independent as 
$b$ varies (notice that ${\cal H}(\TT)$ is countable). For
all $t\geq 0$ define an equivalence relation $\sim_t$ on $\TT$ by 
saying that $v\sim_t w$ if and only if the path $[[v,w]]$ does not contain
any hub $b$ for which $e_b<tL(b)$.
Alternatively, following 
more closely the spirit of Aldous-Pitman's fragmentation, we can 
also say that we consider Poisson point process $(b(t),t\geq 0)$ on the set of
hubs with intensity 
$\d t\otimes\sum_{b\in{\cal H}(\TT)}L(b)\delta_{b}(\d v)$, and for each $t$ we
let 
$v\sim_t w$ if and only if no atom of the Poisson process that has 
appeared before time $t$ belongs to the path $[[v,w]]$. We let 
$\TT^t_1,\TT^t_2,\ldots$ be the distinct equivalence classes for 
$\sim_t$, ranked according to the decreasing order of their $\mu$-masses
(provided these are well-defined quantities). 
It is easy to see that these sets are trees (in the same sense as $\TT$),
and that the families $(\TT^t_i,i\geq 1)$ are nested as $t$ varies, that is, 
for every $t'>t$ and $i\geq 1$, there exists $j\geq 1$ such that 
$\TT^{t'}_i\subset \TT^t_j$. 
If we let $F^+(t)=(\mu(\TT^t_1),\mu(\TT^t_2),\ldots)$, $F^+$ is thus a 
fragmentation process in the sense that $F^+(t')$ is obtained by splitting 
at random the elements of $F^+(t)$.
We mention that the fragmentation $F^+$ is also considered and studied in 
the work in preparation \cite{abdel}, with independent methods. 

We now state our main result, postponing definitions and 
properties of stable subordinators to the next section. Let
$$D_{\alpha}=\frac{\alpha(\alpha-1)\Gamma\left(1-\frac{1}{\alpha}\right)
}{\Gamma(2-\alpha)}
=\frac{\alpha^2\Gamma\left(2-\frac{1}{\alpha}\right)}{\Gamma(2-\alpha)}.$$

\begin{thm}\label{T1}
The process $F^+$ is a self-similar fragmentation with index 
$1/\alpha\in(1/2,1)$
and erosion coefficient $c=0$. Its dislocation measure $\nu_{\alpha}$ is 
characterized by 
$$\nu_{\alpha}(G)=D_{\alpha}E\left[T_1 G(T_1^{-1}\Delta T_{[0,1]})\right]$$
for any positive measurable function $G$, where $(T_x,0\leq x\leq 1)$ 
is a stable subordinator with index $1/\alpha$, characterized by the Laplace
transform 
$$E[\exp(-\lambda T_1)]=\exp(-\lambda^{1/\alpha})\quad \quad \lambda\geq 0,$$
and $\Delta T_{[0,1]}$ is the sequence of the jumps of $T$, ranked by 
decreasing order of magnitude.  
\end{thm}

In a companion 
paper \cite{mierfmoins}, we studied a self-similar 
fragmentation process $(F^-(t),t\geq 0)$ 
which consisted in the decreasing sequences of the $\mu$-masses of 
the connected components of the set $\{v\in\TT:\HT(v)>t\}$ at time $t$, i.e.\ 
the forest obtained by putting aside the vertices of the stable tree  
height less than $t$. 
This fragmentation was studied in the Brownian case by 
Bertoin \cite{bertsfrag02}, although this work does not mention trees and only
uses the encoding height process, which is well-known to be twice the standard
Brownian excursion, and it was showed that it was self-similar with
characteristics ($-1/2,0,\nu_{\rm AP}$) (in \cite{bertsfrag02} the dislocation
measure is found to be $2\nu_{\rm AP}$, but it is done with a different 
normalization, using the standard excursion instead of twice this excursion). 
In \cite{mierfmoins}, we showed that $F^-$ has characteristics
($1/\alpha-1,0,\nu_{\alpha}$), with $\nu_{\alpha}$ as in Theorem
\ref{T1}.  Bertoin's observation that the two devices described above 
for fragmenting the Brownian CRT are ``dual'' (same dislocation measure but
indices with  different signs) is therefore quite  surprisingly generalized in
the larger context of stable trees.  Heuristically, this is made possible by an
exchangeability property of the root of the stable tree with other vertices
(with respect to the measure  $\mu$), which indeed suggests that when removing
a hub or removing the vertices  below a given hub, the subsequent forests will
have the same law up to rescaling. 

Let us now present a second motivation for studying the fragmentation $F^+$. 
As the rest of the paper will show, our proofs involve a lot the theory
of L\'evy processes, and compared with the study of $F^-$, which made a
consequent place to combinatoric tree structures, the study of $F^+$ will 
be mainly ``analytic''. The fact that L\'evy processes may be involved in 
fragmentation processes is not new. According to \cite{bertfrag99} and 
\cite{mier01}, adding a drift to a certain class of L\'evy processes allows
to construct interesting fragmentations related to the entrance boundary of 
the stochastic additive coalescent. Here, rather than adding a drift, which 
by analogy between \cite{jpda97ebac} and \cite{bertfrag99} amounts to 
cut the skeleton of a continuum random tree with a homogeneous Poisson 
process, we will perform a ``removing the jumps'' operation
analog to our inhomogeneous cutting on the hubs of the tree. 

Precisely, let $(X_s,s\geq 0)$ be the canonical process in the 
Skorokhod space $\D([0,\infty),\R)$ and let $P$ be the law of the
stable L\'evy process with index $\alpha\in(1,2)$, upward jumps
only, characterized by the Laplace exponent 
$$E[\exp(-\lambda X_1)]=\exp(\lambda^{\alpha}).$$
As we will recall from the work of Chaumont \cite{chaumont97}
in the following section, we may define the law $N^{(1)}$ of 
the excursion with unit duration
of this process above its infimum process. Under this law, $X_s=0$ for 
$s>1$, so we let $\Delta X_{[0,1]}$ be the sequence of the jumps
$\Delta X_s=X_s-X_{s-}$ for $s\in(0,1]$, ranked in decreasing order of 
magnitude.  
Consider the following marking process on the jumps: conditionally on $X$, let 
$(e_s, s: \Delta X_s>0)$ be a family of 
independent random variables with standard 
exponential distribution, indexed by the countable set of jump-times of $X$. 
For every $t\geq 0$ let
$$Z^{(t)}_s=\sum_{0\leq u\leq s}\Delta X_u \ind_{\{e_u<t \Delta X_u \}}.$$
That is, each jump with magnitude $\Delta$ is marked with 
probability $1-\exp(-t\Delta)$ independently of the other jumps and 
consistently as $t$ varies, and $Z^{(t)}$ is 
the process that sums the marked jumps. We will see that $Z^{(t)}$ is finite
a.s., so we may define $X^{(t)}=X-Z^{(t)}$ under 
$N^{(1)}$. Let 
$$\un{X}^{(t)}_s=\inf_{0\leq u\leq s}X^{(t)}_u \quad,\quad 0\leq s\leq 1,$$
and let $F^{\natural}(t)$ be the sequence of lengths of the constancy 
intervals of the process $\un{X}^{(t)}$, ranked in decreasing order. 

\begin{thm}\label{T2}
The process $(F^{\natural}(t),t\geq 0)$ has the same law as 
$(F^+(t),t\geq 0)$. 
\end{thm}

We organize the paper as follows. 
In Sect.\ \ref{levyproc} we recall some facts about L\'evy
processes, excursions, and conditioned subordinators that will
be  crucial for our study. In Sect.\ \ref{levytree} we 
give the rigorous description of Duquesne and Le Gall's L\'evy trees,
and rephrase the definition of $F^+$ given above in terms of a partition of
the  unit interval associated to a certain marked excursion of a stable
L\'evy  process. Sections \ref{fplus} and \ref{levyrep} are then
respectively dedicated to the study of
$F^+$ and $F^{\natural}$. Asymptotic results are finally given
concerning the behavior at small and large times of
$F^+$ in Sect.\ \ref{smasym}.

\section{Some facts about L\'evy processes}\label{levyproc}

\subsection{Stable processes, inverse subordinators}

Let $(X_s,s\geq 0)$ be the canonical process in the Skorokhod space 
$\D([0,\infty),\R)$ of c\`adl\`ag paths on $[0,\infty)$. We fix
$\alpha\in(1,2)$. Let $P$ be the law on $\D([0,\infty),\R)$ 
that makes $X$ the spectrally positive
stable process with index $\alpha$, that is, $X$ has independent and
stationary increments under $P$, it has only positive jumps, and its
marginal law at some (and then all) $s>0$ has Laplace transform given by
the L\'evy-Khintchine formula:
\begin{equation}\label{laplaX}
E[e^{-\lambda X_s}]=\exp(s \lambda^{\alpha})=
\exp\left(s\int_0^{\infty}\frac{C_{\alpha}\d x}{x^{1+\alpha}}(e^{-\lambda
x}-1+\lambda x)\right),\quad
\lambda\geq 0,
\end{equation} 
where $C_{\alpha}=\alpha(\alpha-1)/\Gamma(2-\alpha)$. 
A fundamental property of $X$ under $P$ is the
{\em scaling property}
$$\left(\frac{1}{\lambda^{1/\alpha}}X_{\lambda s},s\geq 0\right)
\build=_{}^{d}(X_{s},s\geq 0) \quad \mbox{ for all } \lambda>0.$$
We let $(p_s(x),s> 0,x\in\R)$ be the density with respect to Lebesgue
measure of the law $P(X_s\in \d x)$, which is known to exist and to be
jointly continuous in $s$ and $x$. 

Denote by $\un{X}$ the infimum process of $X$ defined by 
$$\un{X}_s=\inf_{0\leq u\leq s}X_u\, ,\quad s\geq 0.$$
Let $T$ be the right-continuous inverse of the increasing process
$-\un{X}$ defined by 
$$T_x=\inf\{s\geq 0: \un{X}_s<-x\}.$$ 
Then it is 
known that under $P$, $T$ is a stable subordinator with index $1/\alpha$, that
is, an increasing L\'evy process with Laplace exponent 
$$E[e^{-\lambda T_x}]=\exp(-x
\lambda^{1/\alpha})=\exp\left(-x\int_0^{\infty}\frac{c_{\alpha}\d
y}{y^{1+1/\alpha}}(1-e^{-\lambda y})\right)\quad\mbox{ for } 
\quad \lambda,x\geq 0, $$
where $c_{\alpha}=(\alpha\Gamma(1-1/\alpha))^{-1}$. 
We denote by $(q_x(s),x,s>0)$ the family of densities with respect to
Lebesgue measure of the law $P(T_x\in\d s)$, by \cite[Corollary
VII.1.3]{bertlev96} they are given by 
\begin{equation}\label{ballotc}
q_x(s)=\frac{x}{s}p_s(-x). 
\end{equation}
We also introduce the notations $P^s$ for the law of of the processes $X$
under $P$, killed at time $s$, and $P^{(-x,\infty)}:=P^{T_x}$ for the 
law of the process killed when it first hits $-x$. 

Let us now discuss the conditioned forms of distributions of jumps of 
subordinators. An easy way to obtain regular versions for 
these conditional laws is developed in \cite{ppy92,pitmanpk02}. First, we
define the size-biased permutation of the sequence $\Delta T_{[0,x]}$ of the
ranked jumps of $T$ in the interval $[0,x]$ as follows. Write
$\Delta T_{[0,x]}=(\Delta_1(x),\Delta_2(x),\ldots)$ with $\Delta_1(x)\geq
\Delta_2(x)\geq\ldots$, and recall that $T_x=\sum_i\Delta_i(x)$. We define, 
following \cite{ppy92,pitmanpk02}, the size-biased ordered sequence
$\Delta_k^*(x),k\geq 1$ as follows. 
Let $1^*$ be a r.v. such that 
$$P(1^*=i| \Delta T_{[0,x]})=\frac{\Delta_i(x)}{T_x}$$ 
for all $i\geq 1$, and set $\Delta^*_1(x)=\Delta_{1^*}(x)$. Recursively, let
$k^*$ be such that
$$P(k^*=i|\Delta T_{[0,x]}, (j^*,1\leq j\leq
k-1))=\frac{\Delta_i(x)}{T_x-\Delta^*_1(x)-\ldots-\Delta^*_{k-1}(x)}$$
for $i\geq 1$ distinct of the $j^*$, $1\leq j\leq k-1$, and finally set
$\Delta^*_k(x)=\Delta_{k^*}(x)$. Then
\begin{lmm}\label{permanform}
{\rm (i)} For $k\geq 1$, 
$$P\left(\Delta^*_k(x)\in \d y\left|
T_x,(\Delta^*_j(x),1\leq j\leq k-1)\right.\right)=\frac{c_{\alpha}xq_x(s-y)}{
sy^{1/\alpha}q_x(s)}\d y$$
where $s=T_x-\Delta^*_1(x)-\ldots-\Delta^*_{k-1}(x)$.  

{\rm (ii)} Consequently, given $T_x=t,\Delta^*_1(x)=y$, the sequence 
$(\Delta^*_2(x),\Delta^*_3(x),\ldots)$ has the same law as $(\Delta^*_1(x),
\Delta^*_2(x)\ldots)$ given $T_x=t-y$. Conversely, if we are given a random 
variable $Y$ with same law as $\Delta^*_1(x)$ given $T_x=t$ and, given 
$Y=y$, a sequence $(Y_1,Y_2,\ldots)$ with same law as 
$(\Delta^*_1(x),\Delta^*_2(x))$ given $T_x=t-y$, then $(Y,Y_1,Y_2,\ldots)$ has 
same law as $(\Delta^*_1(x),\Delta^*_2(x),\ldots)$ given $T_x=t$.  
\end{lmm}

This gives a regular conditional version for $(\Delta^*_i(x),i\geq 1)$
given $T_x$, and thus induces a conditional version for 
$\Delta T_{[0,x]}$ given $T_x$ by ranking. 

\subsection{Marked processes}

We are now going to enlarge the original probability space to 
mark the jumps of the stable process. 
We let $M_X$ be the law of a sequence $\be=(e_s,s:\Delta X_s>0)$ of independent
standard exponential random variables, indexed by the (countable)
set of times where the canonical process $X$ jumps\footnote{One way to attach such 
variables in a measurable way to the $\omega$-dependent set of times $\{s:\Delta X_{s}>0\}$ is to consider 
a doubly-indexed family $(e_{i,j},i,j\geq 1)$ of iid standard exponential variables independent of $X$,
and to attach $e_{i,j}$ to the time of occurrence of the $i$-th largest jump of $X$ in the interval
$[j-1,j)$.}. We let $\bP(\d X,\d\be)=
P(\d X)\otimes M_X(\d\be)$. This probability allows to mark the jumps 
of $X$, precisely we say that a jump occurring at time 
$s$ is marked at level $t\geq 0$ if $e_s<t\Delta X_s$. Write
$$Z^{(t)}_s=\sum_{0\leq u\leq s}\Delta X_u \ind_{\{e_u<t 
\Delta X_u)\}}$$
for the cumulative process of marked jumps at level $t$. We also 
let $X^{(t)}=X-Z^{(t)}$. 
We know that
the process $(\Delta X_s,s\geq 0)$ of the jumps of $X$ is under $P$ 
a Poisson point process with intensity 
$C_{\alpha}x^{-1-\alpha}\d x$ on $(0,\infty)$, it is then standard 
that the process $(\Delta Z^{(t)}_s,s\geq 0)$ is a Poisson point process 
with intensity $C_{\alpha}x^{-\alpha-1}(1-e^{-tx})\d x$, meaning that
under $\bP$, 
$Z^{(t)}$ is a subordinator with no drift and L\'evy measure
$C_{\alpha}x^{-\alpha-1}
(1-e^{-tx})\d x$, more precisely its Laplace transforms 
are given by 
$$\bE[e^{-\lambda Z^{(t)}_s}]=\exp\left(-s
\int_0^{\infty}C_{\alpha}(1-e^{-tx})
\frac{1-e^{-\lambda x}}{x^{\alpha+1}}\d x\right)=\exp(-s(\lambda+t)^{\alpha}
+s\lambda^{\alpha}+st^{\alpha}).$$
We denote by $(\rho^{(t)}_s(x),s,x\geq 0)$ the densities of the laws 
$P(Z^{(t)}_s\in \d x)$. It can be checked by \cite[Proposition 28.3]{Sato} 
from the expression of the L\'evy measure of $Z^{(t)}$ 
that these densities exist and are jointly continuous. 
Likewise, the process $X^{(t)}$ is under $\bP$ a L\'evy process
with L\'evy measure $C_{\alpha}e^{-tx}x^{-\alpha-1}\d x$, 
and the Laplace transform of $X^{(t)}_s$ is given by
$$\bE[e^{-\lambda
X^{(t)}_s}]=\exp\left(s\lambda \alpha t^{\alpha-1}
+s\int_0^{\infty}C_{\alpha}e^{-tx}\frac{\d x}{x^{\alpha+1}}
(e^{-\lambda x}-1+\lambda x)\right)=\exp(s(\lambda+t)^{\alpha}-st^{\alpha}), 
$$ 
which is obtained by dividing the Laplace exponent of $X_s$ by that
of $Z^{(t)}_s$. 

We now state an absolute continuity result that is analogous to 
Cameron-Martin's formula for Brownian motion with drift. 

\begin{prp}\label{abscont}
For every $t,s\geq 0$, we have the following absolute continuity relation:
for every positive measurable functional $F$, 
$$\bE[F(X^{(t)}_u,0\leq u\leq s)]=E[\exp(-st^{\alpha}-tX_s)F(X_u,0\leq u\leq 
s)].$$
\end{prp}

\proof
By the expression for the Laplace exponent of $X^{(t)}$, we get
$$\bE[e^{-\lambda X^{(t)}_s}]=e^{-st^{\alpha}}E[e^{-(\lambda+t)X_s}],$$
hence giving $\bP(X^{(t)}_s\in \d x)=e^{-st^{\alpha}-tx}P(X_s\in \d x)$. The 
result easily follows by the Markov property. \cq

\rem
Such  an  identity   is  a  special  case  of  the   so  called  {\em  density
transformations} for L\'evy processes, see e.g.\ 
\cite[Theorem 33.2]{Sato}.

As a first consequence, it immediately follows that $X^{(t)}$ also has
jointly continuous densities under $\bP$, which are given by 
$$p^{(t)}_s(x)=\frac{\bP(X^{(t)}_s\in \d x)}{\d x}=\exp(-st^{\alpha}-tx)
p_s(x).$$

We let $\un{X}^{(t)}$ be the infimum process of $X^{(t)}$ and 
$T^{(t)}$ the right-inverse process of $-\un{X}^{(t)}$, defined as we did
above define $\un{X}$ and $T$. 

It is easily obtained that for every $t\geq 0$, the process 
$(X,Z^{(t)})$ is again a L\'evy process under the law $\bP$. 
%
%
We will also denote by $\bP^s$, $\bP^{(-x,\infty)}$ the laws derived from 
$P^s$ and $P^{(-x,\infty)}$ by marking the jumps with $M_X$; 
$Z^{(t)}$ and $X^{(t)}$ are then defined as before. 

\subsection{Bridges, excursions}
For $r\in\R$ and $s>0$ we will denote by $P_{0\to r}^s$ 
the law of the {\em stable bridge from $0$ to $r$ with length $s$}, 
so the family $(P_{0\to r}^s,r\in\R)$ forms a regular conditional version for 
$P^s(\cdot| X_s=r)$. 
By \cite{fpy92}, a regular version (which is the one we will always 
consider) is obtained
as the unique law on the Skorokhod space $\D([0,s],\R)$ that 
satisfies the following absolute continuity relation: for every 
$a\in(0,s)$ and any continuous functional $F$,
\begin{equation}\label{absbridge}
P_{0\to r}^s(F(X_u,0\leq u\leq s-a))=E\left[F(X_u,0\leq u\leq s-a)\frac{p_a(
r-X_{s-a})}{p_s(r)}\right].
\end{equation}
We 
let $\bP^s_{0\to r}$ be the marked analog of $P^s_{0\to r}$ on an enriched 
probability space. Notice that Proposition 
\ref{abscont} immediately implies that the bridge laws for 
the process $X^{(t)}$ under $\bP$ are the same as those of $X$. 
Stable bridges from $0$ to $0$ satisfy the following scaling property: under 
$P^v_{0\to 0}$, the process $(v^{-1/\alpha}
X_{vs},0\leq s\leq 1)$ has law $P^1_{0\to 0}$. 

\begin{lmm}\label{markbr}
The following formula holds for any positive measurable $f,g,H$:
\begin{eqnarray*}
\lefteqn{E^1_{0\to0}\left[H(X)\sum_{0\leq s\leq 1}\!\Delta X_s
f(s)g(\Delta X_s)\right]}\\
&=&\int_0^1\!\! \d s\, f(s)\int_0^{\infty}\!\! \d x
\frac{C_{\alpha}p_1(-x)}{x^{\alpha}p_1(0)}g(x)
E^1_{0\to -x}[H(X\oplus(s,x))],
\end{eqnarray*}
where $X\oplus (s,x)$ is the process $X$ to which has been added a jump at
time 
$s$ with magnitude $x$. Otherwise said, a stable bridge from $0$ to $0$ 
together with a jump
$(s,\Delta X_s)$ picked according to the $\sigma$-finite measure
$m(\d s,\d x)=\sum_{u:\Delta X_u>0}\Delta X_u \delta_{(u,\Delta X_u)}
(\d s,\d x)$
is obtained by taking a stable bridge from $0$ to $-x$ and adding a jump with magnitude $x$ at time $s$, where $s$ is uniform in $(0,1)$ and 
$x$ is independent with $\sigma$-finite ``law'' $C_{\alpha}p_1(-x)p_1(0)^{-1}
x^{-\alpha}\d x$.
\end{lmm}

\proof
By the L\'evy-It\^o decomposition of L\'evy processes, one can write, under 
$P$, that $X_s$ is the compensated sum 
$$X_s=\lim_{\eps\to0}\left(\sum_{0\leq u\leq s}\!\Delta X_u 
\ind_{\{\Delta X_u>\eps\}}-(\alpha-1)^{-1}
C_{\alpha}\eps^{1-\alpha} s\right),\quad s\geq 0,$$
where $(\Delta X_u,u\geq 0)$ is a Poisson point process with intensity 
$C_{\alpha}x^{-\alpha-1}\d x$, and where the convergence is almost sure. 
By the Palm formula for 
Poisson processes, we obtain that for positive measurable $f,g,h,H$:
\begin{eqnarray*}
\lefteqn{E^1\left[h(X_1)H(X)\sum_{0\leq s\leq 1}\!\Delta X_s
f(s)g(\Delta X_s)\right]} \\
&=&\int_0^1\!\! \d s\, f(s)\int_0^{\infty}\!\!\d x
\frac{C_{\alpha}}{x^{\alpha}}g(x)
E^1[h(x+X_1)H(X\oplus(s,x))]. 
\end{eqnarray*}
The result is then obtained by disintegrating with respect to the law of 
$X_1$. 
\cq

We now state a useful decomposition of the stable bridge from $0$ to $0$. 
Recall that $(\rho^{(t)}_s(x),x\geq 0)$ is the density of $Z^{(t)}_s$ under 
$\bP$ and that $X^{(t)}_1+Z^{(t)}_1=X_1$, which is a sum of two independent
variables. From this we conclude that 
$(p_1(0)^{-1}p_1^{(t)}(x)\rho^{(t)}_1(x),x\geq 0)$ is a probability density on 
$\R_+$. 

\begin{lmm}\label{decompX}
Take a random 
variable ${\cal Z}$ with law $P({\cal Z}\in \d z)=p_1^{(t)}(-z)
\rho^{(t)}_1(z)p_1(0)^{-1}\d z$. Conditionally on ${\cal Z}=z$, take 
$X'$ with law $P^1_{0\to -z}$ and $Z$ with law
$\bP^1(Z^{(t)}\in \cdot| Z^{(t)}_1=z)$, independently. 
That is, $Z$ is the bridge of 
$Z^{(t)}$ with length $1$ from $0$ to $z$. 
Then $X'+Z$ has law $P^1_{0\to0}$.
\end{lmm}

\rem
The definition for the bridges of $Z^{(t)}$ under $\bP^1$ has not been
given  before. One can either follow an analogous definition as
(\ref{absbridge}),  or use Lemma \ref{permanform} about conditioned jumps
of subordinators.  We explain  this for bridges  of $T$, the  construction for
bridges of $Z^{(t)}$ being similar. Take $(\Delta_i,i\geq 1)$ a sequence whose
law is that of the  jumps $\Delta T_{[0,1]}$ of $T$ under $P$ before time
$1$, ranked in decreasing order, and conditioned by $T_1=z$, in the sense
of Lemma \ref{permanform}. Take also a sequence $(U_i,i\geq 1)$ of
independent uniformly distributed  random variables on $[0,1]$,
independent of $\Delta T_{[0,1]}$.  Then one checks from the L\'evy-It\^o
decomposition for L\'evy processes  that the law $Q_z$ of the process
$T^{\rm br}_s=\sum \Delta_i \ind_{\{s\geq U_i\}}$, with $0\leq s\leq 1$,
defines as $z$ varies a regular version of the  conditional law
$P^1(T\in \cdot| T_1=z)$. 

\proof
Recall that under $\bP^1$, $X$ can be written as $X^{(t)}+Z^{(t)}$ with 
$X^{(t)}$ and $Z^{(t)}$ independent. Consequently, for $f$ and $G$ 
positive continuous, we have 
$$
E^1[f(X_1)G(X)]=\bE^1[f(X^{(t)}_1+Z^{(t)}_1)G(X^{(t)}+Z^{(t)})]$$
so 
\begin{eqnarray*}
\int_{\R}\d x\, p_1(x)f(x)E^1_{0\to x}[G(X)]  
&=&\int_{\R}
\d x \, p_1(x)\int_0^{\infty}\d z \frac{p^{(t)}_1(x-z)\rho^{(t)}_1(z)}{p_1(x)}
f(z)\\
& & \times\bE^1[G(X^{(t)}+Z^{(t)})|X^{(t)}_1=x-z, Z^{(t)}_1=z].
\end{eqnarray*}
Thus, for (Lebesgue) almost every $x$, the bridge with law $P^1_{0\to x}$
is obtained by taking a bridge of $X^{(t)}$ (or $X$ by previous remarks)
from $0$ to $-{\cal Z}_x$ and an independent bridge of $Z^{(t)}$ from 
$0$ to ${\cal Z}_x$, where ${\cal Z}_x$ is a r.v.\ with law 
$\d z\, p_1(x)^{-1}p^{(t)}_1(x-z)\rho^{(t)}_1(z)$ on $\R_+$.  We extend this 
result to every $x\in \R$ by an easily checked continuity result 
for the laws of bridges which stems from (\ref{absbridge}) and the 
continuity of the densities. Taking $x=0$ gives the result. \cq

We now turn our attention to excursions. 
The fact that $X$ has no negative jumps implies that 
$-\un{X}$ is a local time at $0$ for the reflected process 
$X-\un{X}$. Let $N$ be the It\^o excursion measure
of $X-\un{X}$ away from $0$, so that the path of $X-\un{X}$ is obtained 
by concatenation of the atoms of a Poisson measure with intensity 
$N(\d X)\otimes \d t$ on $\D^{\dagger}([0,\infty),\R)\times \R_+$, where 
$\D^{\dagger}([0,\infty),\R)$ denotes the Skorokhod space of paths that 
are killed at some time $\zeta$. Under $N$, almost every path 
$X$ starts at $0$, is positive on an interval $(0,\zeta)$ and dies 
at the first time $\zeta(X)\in(0,\infty)$ 
it hits $0$ again. We let $\bN$ be the enriched law
with marked jumps. It follows from excursion theory that the
L\'evy process $(X,Z^{(t)})$ under
$\bP$ is obtained by taking a Poisson point measure 
$\sum_{i\in I}\delta_{X^i,\be^i,s^i}$ indexed by a countable set $I$, 
with intensity $\bN(\d X,\d \be)
\otimes\d s$, writing $Z^{(t),i}$ for the cumulative process of marked jumps 
for $X^i$ and letting 
$$X_s=-s^i+X^i\left(s-\sum_{j:s^j<s^i} \zeta_j(X^j)\right)$$ 
and 
$$Z^{(t)}_s=\sum_{j:s^j<s^i}Z^{(t),j}_{\zeta_j(X^j)}+Z^{(t),i}\left(s-
\sum_{j:s^j<s^i} \zeta_j(X^j)\right),$$
whenever $\sum_{j:s^j<s^i}\zeta_j(X^j)\leq s\leq 
\sum_{j:s^j\leq s^i}\zeta_j(X^j)$.

If $X$ is stopped at some time $s$, for any $u\in [0,s]$ we define 
the rotated process 
$$V_u X(r)=(X_{r+u}-X_u)\ind_{\{0\leq u< s-u\}}
+(X(r-s+u)+X_s-X_u)\ind_{\{s-u\leq r\leq s\}}.$$
Let $m_s=-\un{X}_s$ and suppose that this minimum is attained only once on 
$[0,s]$. We define the Vervaat transform of $X$ as 
$VX =V_{T(m_s-)}X$, the rotation of $X$ at the time where it attains its 
infimum. Provided that 
$X_0=0$ and $X_s=X_{s-}=0$ (say that $X$ is a bridge), $VX$ is then an 
excursion-like function, starting and ending at $0$, and staying positive 
in the meanwhile.  

We will denote by $N^{(v)}$ the law of $VX$ under $P_{0\to0}^v$, and 
$\bN^{(v)}$ the corresponding ``marked'' version. Call it the 
law of the {\em excursion of $X$ with duration $v$}. The ``Vervaat theorem'' 
in \cite{chaumont97} shows that 
$N^{(v)}$ is indeed a regular conditional version 
for the ``law'' $N(\cdot|\zeta=v)$: for any positive measurable functional
$F$ and function $f$, 
$$N(F(X_s,0\leq s\leq \zeta)f(\zeta))=\int_{\R_+}f(\zeta)N(\zeta\in \d v)
N^{(v)}(F(X_s,0\leq s\leq v)).$$
As for bridges, we 
also have the scaling property at the level of 
conditioned excursions: under $N^{(v)}$, 
$\left(v^{-1/\alpha}X_{v s},0\leq s\leq 1\right)$
has law $N^{(1)}$.
Notice also (either by Vervaat's theorem or 
directly, using Proposition \ref{abscont}) that the excursions
of $X^{(t)}$ under $\bP$, conditioned to have a fixed duration $v$ 
are the same as that of $X$ under $N^{(v)}$.

\section{The stable tree}\label{levytree}

\subsection{Height Process, width process}

We now introduce the rigorous definition and useful properties of 
the stable tree.  
This section is mainly inspired by \cite{duqleg02,duq02}. With the notations 
of section \ref{levyproc}, and for 
$t\geq 0$, let $R^{(t)}$ be the time-reversed process of $X$ at time $t$:
$$R^{(t)}_s=X_t-X_{(t-s)-}\quad 0\leq s\leq t.$$
It is standard that this process has the same law as $X$ killed
at time $t$ under $P$. Let $\ov{R}^{(t)}$ be its supremum process, and 
$\wh{L}^{(t)}$ be the local time process at level $0$ of the 
reflected process $\ov{R}^{(t)}-R^{(t)}$. 
We let $H_t=\wh{L}^{(t)}_t$. The normalization for $\wh{L}^{(t)}$
is chosen so that 
$$H_t=\lim_{\eps\downarrow 0}\frac{1}{\eps}\int_0^t\ind\{
\ov{R}^{(t)}_s-R^{(t)}_s\leq \eps\}\d s,$$
in probability for every $t$. It is proved in \cite{duqleg02} that 
$H$ admits a continuous modification, which is the one we are going 
to work with from now on. 
It has to be noticed that $H_t$ is not a Markov process, except in the 
case where $X$ is Brownian motion. As a matter of fact, it can be noticed 
that $H$ admits infinitely many local minima attaining the same value
as soon as $X$ has jumps. To see this, consider a jump time $t$ of $X$, and 
let $t_1,t_2>t$ so that $\inf_{t\leq u\leq t_i}X_u=X_{t_i}$ and 
$X_{t-}<X_{t_i}<X_t$, $i\in\{1,2\}$. Then it is easy to see that 
$H_{t}=H_{t_1}=H_{t_2}$ and that one may in fact find an infinite number
of distinct $t_i$'s satisfying the properties of $t_1,t_2$. On the other hand, 
it is not difficult to see that $H_t$ is a local minimum of $H$, see
Proposition \ref{onetoone} below. 

It is  shown in \cite{duqleg02} that the definition of $H$ 
still makes sense under the 
$\sigma$-finite measure $N$ rather than the probability law $P$. The process
$H$ is then defined only on $[0,\zeta]$, and we call it the {\em excursion
of the height process}.
Using the scaling property, one can then define
the height process under the laws $N^{(v)}$. Call it the law of the 
excursion of the height process with duration $v$.

The key tool for defining the local time of hubs is the 
local time process of the height process. We will denote by 
$(L^t_s,t,s\geq0)$. It can be obtained a.s.\ for every fixed $s,t$ by 
\begin{equation}\label{loctimeapp}
L^t_s=\lim_{\eps\to0}
\frac{1}{\eps}\int_0^s\ind_{\{t< H_u\leq t+\eps\}}\d u.
\end{equation}
That is, $L_s^t$ is the density of the occupation measure of $H$ at level
$t$ and time $s$. For $t=0$,
one gets $(L_s^0,s\geq 0)=(\un{X}_s,s\geq 0)$, 
which is a reminiscent of the fact that the excursions of the height 
process are in one-to-one correspondence with excursions of $X$ with the 
same lengths.

It is again possible to define
the local time process under the excursion measures $N$ and $N^{(v)}$. 
Duquesne and Le Gall \cite{duqleg02} have shown that under $P$, the process 
$(L_{T_x}^t,t\geq 0)$ has the law of the continuous-stable branching process
starting at $x>0$, with stable ($\alpha$)
branching mechanism. 
One can get interpretations for the process $(L_{\zeta}^t,
t\geq 0)$ under the measure $N$ or of $(L_v^t,t\geq 0)$ under $N^{(v)}$ in
terms of conditioned continuous-state branching processes, see 
\cite{mierfmoins}. 

\subsection{The tree structure}\label{secsttree}

Let us motivate the term of ``height process'' for $H$ by embedding 
a tree inside $H$, following \cite{legall91,aldouscrt93}. Consider the 
height process $H$ under the law $N^{(1)}$.  
We can define a pseudo metric $D$ on $[0,1]$ by letting
$D(s,s')=H_s+H_{s'}-2\inf_{u\in [s,s']}H_u$ (with the convention 
that $[s,s']=[s',s]$ if $s'<s$). Let $s\equiv s'$ if and only if $D(s,s')=0$. 
\begin{defn}\label{sttree}
The {\em stable tree} $(\TT,d)$ is the quotient of the pseudo-metric space
$([0,1],D)$ by $\equiv$. The {\em root} of $\TT$ is the equivalence class of 
$0$. The {\em mass measure} $\mu$ is the Borel 
measure induced on $\TT$ by Lebesgue's measure 
on $[0,1]$ (so its support is $\TT$). 
\end{defn}
In the sequel, we will often identify $\TT$ with $[0,1]$, even if the 
correspondence is not one-to-one. 
Some comments on this definition. First, the way the tree is embedded in the 
function $H$ can seem quite intricate. It is not difficult, however, to 
see what its ``marginals'' look like. For any finite set of vertices
$s_1,s_2,\ldots,s_k\in [0,1]$, one recovers the structure of the subtree 
spanned by the root and $s_1,s_2,\ldots,s_k$, according to the following 
simple rules: 
\begin{itemize}
\item The height of $s$ is $\HT(s)=H_{s}$.
\item The common ancestor of $s_1,\ldots,s_k$ is 
$b=b(s_1,\ldots,s_k)\in[\min_{1\leq i\leq k}s_i,\max_{1\leq i\leq k}s_i]$ 
such that $H_b=\inf\{H_s:s\in 
[\min_{1\leq i\leq k}s_i,\max_{1\leq i\leq k}s_i]\}$.
\end{itemize}
Notice that all such $b$ are equivalent with respect to $\equiv$. 
The fact that $(\TT,d)$ is indeed a tree (a complete metric space such that 
the only simple path leading from a vertex to another is the geodesic) is intuitive 
and proven in \cite{duqlegprep}. It follows from the construction of 
``marginals'' of $\TT$ in \cite{duqleg02} that given $\mu$, 
$\mu$-a.e.\ vertex is a leaf of $\TT$.


We now relate properties on the stable tree to path properties of the 
underlying L\'evy process we started with to construct the 
height process. We understand here that $X$ and $H$ are defined under 
$N^{(1)}$. Recall that $\TT_b$ stands for the fringe subtree 
rooted at $b$. 

\begin{prp}\label{onetoone}
{\rm (i)} Each hub $b\in{\cal H}(\TT)$ 
is encoded by exactly one time $\tau(b)\in[0,1]$ such that 
$L(b)=\Delta X_{\tau(b)}>0$, and $L(b)$ is given by (\ref{widthhub}) a.s.

{\rm (ii)} If 
$\sigma(b)=\inf\{s\geq \tau(b): X_s=X_{\tau(b)-}\}$,
then 
$\TT_b=[\tau(b),\sigma(b)]/\equiv$. 

{\rm (iii)} 
More precisely, let $\TT^1_b,\TT^2_b,\ldots$ be the connected components 
of $\TT_b\setminus\{b\}$, arranged in decreasing order of mass.  
Let $([\tau_i(b),\sigma_i(b)],i\geq 1)$ 
be the constancy intervals of the infimum 
process of $(X_{s}-X_{\tau(b)},\tau(b)\leq s\leq \sigma(b))$, 
and ranked in decreasing order of length.  
Then $\TT^i_b=(\tau_i(b),\sigma_i(b))/\equiv$. 
\end{prp}

\proof
(i)
Working first 
under $P$, fix $\ell>0$ and let $\tau_{\ell}=\inf\{s\geq 0:\Delta X_s>\ell\}$.
Then $\tau_{\ell}$ is a stopping time for the natural filtration 
associated to $X$, as well as $\sigma_{\ell}=\inf\{s>\tau_{\ell}: 
X_s=X_{\tau_{\ell-}}\}$. By the Markov property, the process 
$X_{[\tau_{\ell},\sigma_{\ell}]}=(X_{\tau_{\ell}+s}-X_{\tau_{\ell}},0\leq
s\leq \sigma_{\ell}-\tau_{\ell})$ is independent of 
$(X_{s+\sigma_{\ell}}-X_{\sigma_{\ell}},s\geq 0)$, which has the same law as 
$X$, and of $(X_s,0\leq s\leq \tau_{\ell})$ conditionally on its final jump 
$\Delta X_{\tau_{\ell}}$. Now if we remove this jump, 
that is, if we let $(\wt{X}_s,0\leq s\leq \tau_{\ell})$ be the modification of 
$(X_s,0\leq s\leq \tau_{\ell})$ that is left-continuous at $\tau_{\ell}$, 
then $\wt{X}$ has the law of a stable L\'evy process killed at some 
independent exponential time, and conditioned to have jumps with 
magnitude less than $\ell$. Also, conditionally on 
$\Delta X_{\tau_{\ell}}=x$, $X_{[\tau_{\ell},\sigma_{\ell}]}$ has the 
law $P^{(-x,\infty)}$ of the stable process killed when it first hits $-x$. 
Hence, by the additivity of the local time and the definition of $H$, 
one has that for every $s\in[\tau_{\ell},\sigma_{\ell}]$, 
$H_s=H_{\tau_{\ell}}+\wt{H}_{s-\tau_{\ell}}$, where $\wt{H}$ is an 
independent copy of $H$, killed when its local time at $0$ attains $x$. 
Consequently, one has $H_s\geq H_{\tau_{\ell}}$ for every 
$s\in[\tau_{\ell},\sigma_{\ell}]$ and $H_{\sigma_{\ell}}=H_{\tau_{\ell}}$, 
moreover, 
one has that for every $\eps>0$, 
\begin{equation}\label{hubrel}
\inf_{(\tau_{\ell}-\eps)\vee 0\leq s\leq \tau_{\ell}}H_s
\vee \inf_{\sigma_{\ell}\leq s\leq \sigma_{\ell}+\eps}H_s<H_{\tau_{\ell}},
\end{equation}
as a consequence of the following fact. By the left-continuity of 
$X$ at $\tau_{\ell}$, for any $\eps>0$ we may find 
$s\in[\tau_{\ell}-\eps,\tau_{\ell}]$ such that 
$\inf_{u\in[s,\tau_{\ell}]}X_u=X_s$. This implies $H_s=H_{\tau_{\ell}}-
\wh{L}^{(\tau_{\ell})}_{
\tau_{\ell}-s}$, and this last term is a.s.\ strictly less than 
$H_{\tau_{\ell}}$ because $0$ is is a.s.\ not a holding point for 
$(\wh{L}^{(\tau_{\ell})}_s, 0\leq s\leq \tau_{\ell})$. This last fact 
is obtained by a 
time-reversal argument, using the fact that the points of increase of the 
local time $\wh{L}^{(t)}$ correspond to that of the supremum process of 
$R^{(t)}$. Moreover, the fact that $X$ has
only positive jumps under $P$ implies that for some suitable $\eps'>0$, 
one can find some 
$s'\in [\sigma_{\ell},\sigma_{\ell}+\eps']$ and some 
$s''\in[\tau_{\ell}-\eps,\tau_{\ell}]$ such that 
$H_u\geq H_{s'}=H_{s''}$ for every $u\in[s',s'']$, and such that 
again $\inf_{u\in[s'',\tau_{\ell}]}X_u=X_{s''}$. Thus the claimed inequality. 
In terms of the structure of the stable tree, (\ref{hubrel}) implies that 
a branchpoint $b$ of the tree is present at height $H_{\tau_{\ell}}$, which is encoded by all the 
$s\in[\tau_{\ell},\sigma_{\ell}]$ such that $H_s=H_{\tau_{\ell}}$, i.e.\ 
such that $X_s=\inf_{u\in [\tau_{\ell},s]}X_u$ (there is always an infinite
number of them).  By definition, the mass measure of the vertices in $\TT_b$ at
distance less than $\eps$ of $b$ is exactly the Lebesgue measure of 
$\{s\in[\tau_{\ell},\sigma_{\ell}]:\wt{H}_{s-\tau_{\ell}}<\eps\}$. Thus by
(\ref{loctimeapp}) we can conclude that $L(b)$ defined at (\ref{widthhub})
exists and equals 
$\wt{L}^0_{\sigma_{\ell}-\tau_{\ell}}=x$ where $\wt{L}$ is the
local time associated to $\wt{H}$. The same argument 
allows to handle the second, third, ... jumps that are $>\ell$. Letting
$\ell\downarrow 0$ 
implies that to any jump of $X$ with magnitude $x$ corresponds 
a hub of the stable tree with local time $x$.  
By excursion 
theory and scaling, the same property holds under 
$N$ and $N^{(1)}$. 

Conversely, suppose that $b$ is a branchpoint in the stable tree. This means that 
there exist        times $s_1<s_2<s_3$ such that $H_{s_1}=H_{s_2}=H_{s_3}$ and 
$H_s\geq H_{s_1}$ for every $s\in[s_1,s_3]$. Let 
$$\tau(b)=\inf\{s\leq s_2:H_s=H_{s_2} \mbox{ and }H_u\geq H_{s_2} \forall
u\in[s,s_2]\}$$
and 
$$\sigma(b)=\sup\{s\geq s_2:H_s=H_{s_2} \mbox{ and }H_u\geq H_{s_2} \forall
u\in[s_2,s]\}$$
(which are not stopping times).
If $\Delta X_{\tau(b)}>0$, we are in the preceding case. Suppose that 
$\Delta X_{\tau(b)}=0$, then by the same arguments as above, $X_s\geq 
X_{\tau(b)}$ for $s\in[\tau(b),\sigma(b)]$, else we could find some 
$s'\in[\tau(b),\sigma(b)]$ such that $H_{s'}<H_{\tau(b)}$. 
Also, the points
$s\in[\tau(b),\sigma(b)]$ such that $H_s=H_{\tau(b)}$ must then satisfy
$X_s=X_{\tau(b)}$ (else there would be a strict increase of the local time of 
the reversed process). This implies that $X_{\tau(b)}$ is a local infimum 
of $X$, attained at $s$. By standard considerations, 
such local infima cannot be attained 
more than three times on the interval $[\tau(b),\sigma(b)]$, 
a.s. But if it was attained exactly three times, 
then the branchpoint would have degree $3$, which is impossible according to the 
analysis of $F^-$ in \cite{mierfmoins}, which implies that 
all hubs of the stable tree have infinite degree.

Assertion (ii) follows easily from this, and (iii)
comes from the fact that the points $u\in[\tau(b),\sigma(b)]$ with 
$H_u=H_{\tau(b)}$ are exactly those points where 
$\inf_{r\in[\tau(b),u]}X_r =X_u$, and the definition of the mass measure on
$\TT$.
\cq

\subsection{A second way to define $F^+$}\label{fplusalt}

We will now give some elementary properties of $F^+$ and rephrase its
definition directly from the excursion of the underlying stable
excursion $X$ rather than the tree itself. First recall that given
$\TT$, we defined $F^+$ through a marking procedure on ${\cal H}(\TT)$ by
taking a Poisson process $(b(t),t\geq 0)$ with intensity $\d t\otimes
\sum_{b\in{\cal H}(\TT)}L(b)\delta_{b}(\d v)$, and by saying  that $b$
is marked at level $t$ if $b\in\{b(s),0\leq s\leq t\}$. By proposition 
\ref{onetoone}, $F^+$ can thus be defined
under the marked law $\bN^{(1)}$. To describe this construction a bit more, we
begin with the following

\begin{lmm}\label{finitehub}
Let $s\in[0,1]$, and write $v(s)$ for the vertex of $\TT$ encoded 
by $s$. Then almost-surely, 
$$\sum_{b\in{\cal H}(\TT)\cap[[{\rm root},v(s)]]}L(b)<\infty.$$
In particular, almost surely, for every hub $b\in{\cal H}(\TT)$ and 
$t\geq 0$, there is at most a finite number of 
hubs marked at level $t$ on the path $[[{\rm root},b]]$. 
\end{lmm}

\proof
Let $s$ be the leftmost time in $[0,1]$ that encodes $v$. It follows
from  Proposition \ref{onetoone} (ii) (and the fact that a.s.\ under
$P$,  every excursion of $R^{(s)}$ below $\ov{R}^{(s)}$ ends by a
jump) that the hubs $b$ in the path $[[{\rm root},v]]$ are all
encoded  by the times $s'<s$ such that $\ov{R}^{(s)}$ jumps at time
$s-s'$. This jump  corresponds to a jump of the reversed process
$R^{(s)}$, whose magnitude
$\Delta R^{(s)}_{s-s'}\geq \Delta \ov{R}^{(s)}_{s-s'}$ equals $L(b)$ by 
Proposition \ref{onetoone} (i). Therefore, we have to show that the
sum of these jumps  is finite a.s. By excursion theory and
time-reversal, it suffices to show that  under $P$, letting $\ov{X}$
be the supremum process of $X$,
\begin{equation}\label{sumfin}
\sum_{0\leq s'\leq s: \Delta\ov{X}_{s'}>0}\Delta X_{s'}<\infty\, ,\qquad
s\geq 0. 
\end{equation}
Now by excursions and Poisson processes theories (see e.g.\ Formula 
(10) in the proof of \cite[Lemma 1.1.2]{duqleg02}), after appropriate 
time-change by the inverse local time at $0$ of the process $X-\ov{X}$, 
the jumps $\Delta X_{s'}$ that achieve new suprema form a 
Poisson point process with intensity $x\times C_{\alpha}x^{-1-\alpha}\d
x$. Since this measure integrates $x$ on a neighborhood of $0$, 
the sum in (\ref{sumfin}) is a.s.\ finite. 

The statement on hubs follows since for any hub $b$ 
encoded by a jump-time $\tau(b)$, 
there is a rational number $r'\in[\tau(b),\sigma(b)]$ 
which encodes some vertex $v$ in the fringe subtree rooted at $b$. Therefore, 
almost-surely, for every $b\in{\cal H}(\TT)$, 
the sum of widths of the hubs on the path $[[\varnothing,b]]$ is finite. It is
then easy to check that if $(x_1,x_2,\ldots)$ is a sequence with finite sum and
if the $i$-th term is marked with probability $1-e^{-t x_i}$, then a.s.\ 
only a finite number of terms are marked. Therefore, 
a.s.\ for every $b\in{\cal H}(\TT)$, there is 
only a finite number of marked hubs on the path
$[[\varnothing,b]]$. \cq

By
definition, two vertices $v,w\in \TT$ satisfy $v\sim_t w$ if and only
if $\{b(s):0\leq s\leq t\}\cap [[v,w]]=\emptyset$. Let ${\cal
H}_t=\{b(s):0\leq s\leq t\}$. For $b\in{\cal H}_t$,
let $\TT^1_b,\TT^1_b,\ldots$ be the connected components of
$\TT_b\setminus\{b\}$ ranked in decreasing order of total mass. We
know that  these trees are encoded by intervals of the
form $(\tau_i(b),\sigma_i(b))$ whose union is
$[\tau(b),\sigma(b)]\setminus\{u:u\equiv b\}$. Define
$$C(t,b,i)=\TT^i_b\setminus \bigcup_{b'\in{\cal H}_t\cap
\TT^i_b}\TT_{b'},$$
so $C(t,b,i)$ is the connected component of the $i$-th largest subtree growing
from $b$ obtained when the hubs marked at level $t$ are deleted.
Plainly, $C(t,b,i)$ is an equivalence class for $\sim_t$ for every
$b\in{\cal H}_t$ and $i\geq1$.  By (iii) in Proposition 
\ref{onetoone}, with obvious notations,  
$$C(t,b,i)\equiv(\tau_i(b),\sigma_i(b))\setminus
\bigcup_{b'\in 
\TT^b_i\cap{\cal H}_t}[\tau(b'),\sigma(b')].$$ 
We also let 
$C(t,\emptyset)$ be the set of vertices whose path to the root does 
not cross any marked hub at level $t$, which is equivalent to 
$[0,1]\setminus \bigcup_{b\in {\cal H}_t} [\tau(b),\sigma(b)]$.
Then $C(t,\emptyset)$ is also an equivalence class for $\sim_t$. 
Intuitively, the classes $C(t,\emptyset)$ and 
$C(t,b,i)$ for $b$ a hub are the equivalence classes for 
$\sim_t$ that have a positive weight. We will see later that the 
rest is a set of leaves of mass zero.

Let us now translate the relation $\sim_t$ 
in terms of the stable excursion $X$ under $\bN^{(1)}$. Let $s,s'\in[0,1]$
encode respectively the vertices $v\neq w\in\TT$. 
Again by Proposition
\ref{onetoone} (ii), the branchpoint $b(v,w)$ of $v$ and $w$ is
encoded by the largest $u$ such that the processes
$(\ov{R}^{(s)}_{s-u+r},0\leq r\leq u)$ and
$(\ov{R}^{(s')}_{s'-u+r},0\leq r\leq u)$ coincide. Let $u(s,s')$ be
the jump-time of $X$ that encodes this branchpoint. Then $v\sim_t w$
if and only if the (left-continuous) processes 
$(\ov{R}^{(s)}_{s-r},u(v,w)\leq r\leq
s)$ and $(\ov{R}^{(s)}_{s'-r},u(v,w)\leq r\leq s')$ never jump at
times when marked jumps at level $t$ for $X$ occur. 

In particular, we may rewrite the equivalence classes
$C(t,b,i)$ and $C(t,\emptyset)$ as follows.  Let $z^t_1\geq z^t_2\geq
\ldots\geq 0$ be the marked  jumps of $X$ at level $t$ under
$\bN^{(1)}$, ranked in decreasing order, and let
$\tau^t_1,\tau^t_2,\ldots$ the corresponding jump times (i.e.\ such
that $\Delta Z^{(t)}_{\tau^t_i}=z^t_i$). For every $i$, let 
$$\sigma^t_i=\inf\{s>\tau^t_i: X_s=X_{\tau^t_i-}=X_{\tau^t_i}-z^t_i\}$$
be the first return time to level $X_{\tau^t_i-}$ after time 
$\tau^t_i$. Define the intervals $I^t_i=[\tau^t_i,\sigma^t_i]$, 
so $I^t_i/\equiv$ is the fringe subtree of the marked hub that has
width $z^t_i$. Notice that the $I^t_i$'s are by no means disjoint, since these
fringe subtrees contain other marked hubs, that might even have greater width. 
For each $i$, the jump with magnitude
$z^t_i$ gives rise to a family of excursions of $X$ above its
minimum. Precisely, let $(X^t_{i,1},X^t_{i,2},\ldots)$ the sequence
of excursions above its infimum of the process 
$$X^t_i(s)=X_{\tau^t_i+s}-X_{\tau^t_i}\quad 0\leq s\leq
\sigma^t_i-\tau^t_i, i\geq 1$$ 
where the $(X^t_{i,j},j\geq 1)$ are arranged by decreasing order of
duration. Let also $I^t_{i,j}=[\tau^t_{i,j},\sigma^t_{i,j}]$  be the
interval in which $X^t_{i,j}$ appears in $X$, so that $\ov{\bigcup_j
I^t_{i,j}}=I^t_i$.  Consider the set
$$C^t_{i,j}=I^t_{i,j}\setminus \bigcup_{k: 
I^t_k\subsetneq I^t_i} I^t_k.$$
By Lemma \ref{finitehub}, there exists some set of indices $k'$ such
that $I^t_{k'}\subsetneq I^t_{i,j}$ and  so that the $I^t_{k'}$'s are
maximal with this property (else we could find  an infinite number of
marked hubs on a path from the root to one of the  hubs encoded by
the left-end of some $I^t_k\subsetneq I^t_{i,j}$).  The Lebesgue
measure of $C^t_{i,j}$ is thus equal to 
$$|C^t_{i,j}|=\sigma^t_{i,j}-\tau^t_{i,j}-\sum
(\sigma^t_k-\tau^t_k),$$ 
where the sum is over the $k$'s such that $I^t_k\subsetneq I^t_i$ and 
the $I^t_k$'s are maximal with this property.
Writing $C^t_0=[0,1]\setminus \bigcup_{i=1}^{\infty}I^t_i$, we
finally get (identifying Borel subsets of $[0,1]$ with Borel subsets
of $\TT$):
\begin{lmm}\label{cij}
The sets $C^t_0$ and $C^t_{i,j}$, for $i,j\geq 1$, are a relabeling of the 
sets $C(t,\emptyset)$ and $C(t,b,i)$.
\end{lmm}

Notice also that another consequence of Lemma
\ref{finitehub} is that $F^+$ is continuous in probability at time
$0$. Indeed, as $t\downarrow 0$, the component $C(t,\emptyset)$ of
the  fragmented tree containing the root increases
to $C(0+,\emptyset)$. Suppose $\mu(C(0+,\emptyset))<1$ with positive
probability. Given $\TT$ take $L_1,L_2,\ldots$ independent  with law
$\mu$. By the law of large numbers, with positive probability a 
positive proportion of the $L_i$'s are separated from the root at
time 
$0+$. However, as a consequence of  
Lemma \ref{finitehub}, a.s.\ for every $n\geq 1$ and $t$ small enough, there 
is no marked hub on the paths $[[{\rm root},L_i]],1\leq i\leq n$, hence a 
contradiction.

\section{Study of $F^+$}\label{fplus}

The goal of this section is to study the fragmentation $F^+$ through the 
representation given in the last section. The first step is to study the 
behavior of the excursion on the equivalence classes $C_{i,j}^t$ and $C_0^t$ 
defined previously. 

\subsection{Self-similarity}

This section is devoted to the proof that $F^+$ is a self-similar fragmentation
with index $1/\alpha$ and no erosion. 

Let us first introduce some notation. Let $(f(x),0\leq x\leq \zeta)\in 
\D^{\dagger}([0,\infty),\R)$ 
be a c\`adl\`ag function with lifetime $\zeta\in[0,\infty)$. By
convention we let $f(x)=f(\zeta)$ for $x>\zeta$. We define the unplugging
operation ${\tt UNPLUG}$ as follows. Let $([a_n,b_n],n\geq 1)$ be a
sequence of disjoint closed intervals with non-empty interior, such that
$0<a_n<b_n<\zeta$ for every $n$. Define the increasing continuous function
$$x^{-1}(s)=s-\sum_{n\geq 1}(s\wedge b_n-a_n)^+\, ,\qquad
s\geq 0,$$
where $a^+=a\vee 0$ and
where the sum converges uniformly on $[0,\zeta]$.
We say that the intervals
$[a_n,b_n]$ are {\em separated} if $x^{-1}(a_n)<x^{-1}(a_m)$ for every $n\neq m$
such that
$a_n<a_m$. This is equivalent to the fact that for every $n\neq m$ with
$a_n<a_m$, the set
$[a_n,a_m]\setminus\bigcup_i[a_i,b_i]$ has positive Lebesgue measure, and it implies 
that the constancy intervals of $x^{-1}$ are exactly $[a_{n},b_{n}],n\geq 1$. If 
$([a_{n},b_{n}],n\geq 1)$ 
is separated, define $x$ as the right-continuous inverse of
$x^{-1}$, then $f\circ x$ is c\`adl\`ag (notice that $(f\circ
x)(s-)=f(x(s-)-)$ for $s\in[0,x^{-1}(\zeta)]$), call it 
${\tt UNPLUG}(f,[a_n,b_n],n\geq 1)$. The
action of ${\tt UNPLUG}$ is thus to remove the bits of the path of $f$ that are
included in $[a_n,b_n]$. 
Last, if we are given intervals $[a_n,b_n]$ that are
not overlapping (i.e.\ such that $a_n<a_m<b_n<b_m$ does not happen for
$n\neq m$, though we might have $[a_n,b_n]\subset[a_m,b_m]$), but such that 
there is a separated subsequence 
$([a_{\phi(n)},b_{\phi(n)}],n\geq 1)$ of maximal intervals 
that covers $\bigcup_n[a_n,b_n]$, we similarly
define the unplugging operation by simply ignoring the non-maximal
intervals. 

\begin{lmm}\label{unpluglmm}
Let $([a_n,b_n],n\geq 1)$ be a sequence of separated intervals, and let 
$\pi$ be a partition of $\N$ with blocks $\pi_{1},\pi_{2},\ldots$. Then, as
$N\to \infty$, ${\tt UNPLUG}(f,[a_n,b_n]:n\in\pi_{1}\cup\ldots\cup\pi_{N})$ 
converges to ${\tt UNPLUG}(f,[a_n,b_n],n\geq 1)$ in the Skorokhod topology. 
\end{lmm}

\proof
Define 
$$x^{-1}_N(s)=s-\sum_{n\in\pi_{1}\cup\ldots\cup\pi_{N}}
(s\wedge b_n-a_n)^+\, ,\qquad
s\geq 0.$$ 
The separation of
intervals ensures that every jump of $x$ corresponds to a jump of $x_N$
for some large $N$, and it is not hard to see that this implies 
$x_{N}(x_{N}^{-1}(x(s)))=x(s)$ for all $s$. Since $f\circ x$ is c\`adl\`ag with duration
$\zeta'=\zeta-\sum_n(b_n-a_n)$, for every $N$ we may find a
sequence of times $0=s_0<s_1<s_2<\ldots<s_{k(N)}=\zeta'$ such that the
oscillation
$$\omega(f\circ x,[s_i,s_{i+1}))=\sup_{s,s'\in[s_i,s_{i+1})} |f\circ x
(s)-f\circ x(s')|\build\to_{N\to\infty}^{}0,$$ 
this uniformly in $1\leq i< k(N)$. 
Let also $s_i^N=x^{-1}_N(x(s_i))$ be the corresponding times
for $f\circ x_N$.
We build a time change
$\lambda_N$ (a strictly increasing continuous function) 
by setting $\lambda_N(s_i)=s^N_i$ for $1\leq i\leq k(N)$, and interpolating
linearly between these times. Easily
$|\lambda_N(s_i)-s_i|\leq\sum_{n\notin\pi_{1}\cup\ldots\cup\pi_{N}}(b_n-a_n)\to 0$, and it follows that
$\lambda_N$ converges pointwise and uniformly to the identity
function of $[0,\zeta']$. On the other hand, $f\circ
x(s_i)=f\circ x_N\circ \lambda_N(s_i)$, so for $s\in(s_i,s_{i+1})$, 
$$|f\circ x_N\circ \lambda_N(s)-f\circ x(s)|\leq \omega(f\circ
x,[s_i,s_{i+1}))+
|f\circ x_{N}\circ\lambda_{N}(s)-f\circ x_{N}\circ \lambda_{N}(s_{i})|.$$
To bound the second term, notice that $x_{N}((s^N_{i},s^N_{i+1}))\subset x((s_{i},s_{i+1}))
\cup\bigcup_{n\notin\pi_{1}\cup\ldots\cup\pi_{N}}[a_{n},b_{n}]$. Therefore 
\begin{eqnarray*}
|f\circ x_{N}\circ\lambda_{N}(s)-f\circ x_{N}\circ \lambda_{N}(s_{i})|&\leq&
\omega(f\circ x,[s_i,s_{i+1}))\\
&+&\sup_{n\notin\pi_{1}\cup\ldots\cup\pi_{N}}
(f(a_n)-f(a_n-)+\omega(f,[a_n,b_n])). 
\end{eqnarray*}
We can conclude that $f\circ x_N\circ\lambda_N$ converges uniformly to
$f\circ x$ since the oscillation $\omega(f,[a_n,b_n])$ converges
to $0$ uniformly in $n\notin \pi_{1}\cup\ldots\cup\pi_{N}$ as $N\to\infty$, 
as does the jump $f(a_n)-f(a_n-)$.
\cq

Under the law 
$\bP^{(-z,\infty)}$ under which $X$ is killed when it first attains $-z$, 
for every $t>0$ we let $z^t_1\geq z^t_2\geq\ldots\geq 0$ be the marked
jumps of $X$ at level $t$, ranked in decreasing order of magnitude, and 
$\tau^t_i$ be the time of occurrence of the jump with magnitude $z^t_i$, while 
$\sigma^t_i$ is the first time after $\tau^t_i$ when $X$ hits 
level $X_{\tau^t_i-}$ (notice that $\tau^t_i, \sigma^t_i$ are not stopping
times). 
Similarly as before, we let $I^t_i=[\tau^t_i,\sigma^t_i]$. 

\begin{lmm}\label{plugg}
For every $z,t\geq 0$, the process ${\tt UNPLUG}(X,(I^t_i : i\geq 1))$
has same law as $X^{(t)}$ under $\bP$, killed
when it first hits $-z$. 
\end{lmm}

Part of this lemma is that it makes sense to apply the unplugging
operation with the intervals $I^t_i$, that is, that these intervals
admit a separated covering maximal sub-family. 

\proof
The fact that the intervals $I^t_i$ admit a covering maximal sub-family 
is obtained by re-using the proof of Lemma \ref{finitehub} and the argument given just
after the definition of $C^t_{i,j}$ in the preceding section. 
Next, write $X=X^{(t)}+Z^{(t)}$. For $a>0$,
let $\tau^{t,a}_1$ be the time of the first jump of $Z^{(t)}$ that is 
$>a$, and let $\sigma^{t,a}_1=\inf\{u\geq \tau^{t,a}_1:X_u=X_{\tau^{t,a}_1-}
\}$. Recursively, let $\tau^{t,a}_{i+1}=\inf\{u\geq \tau^{t,a}_i:
\Delta Z^{(t)}_u>a\}$ and $\sigma^{t,a}_{i+1}=\inf\{u\geq \tau^{t,a}_{i+1}:
X_u=X_{\tau^{t,a}_{i+1}-}\}$. Let
$Z^{(t,a)}_s=\sum_{u\leq s}\Delta Z^{(t)}_u \ind_{\{\Delta Z^{(t)}_u\leq
a\}}$. 
The $\tau^{t,a}_i$'s are stopping times for the filtration generated by
$(X^{(t)},Z^{(t)})$, as well as the
$\sigma^{t,a}_i$'s. By a repeated use of the Markov property at these times 
we get 
$${\tt UNPLUG}(X;(I^t_i:z^t_i>a))\build=_{}^{d}X^{(t)}+Z^{(t,a)},$$ 
where this last process is killed at the time $T^{(t,a)}_z$ when it first
hits $-z$.  In particular, $T_z-\sum_i (\sigma^{t,a}_i-\tau^{t,a}_i)$ has
the same law as $T^{(t,a)}_z$, which converges in law to $T^{(t)}_z$
as $a\downarrow0$ because $Z^{(t,a)}$ converges to
$0$ uniformly on compact sets, and $X^{(t)}$ enters $(-\infty,-z)$ immediately
after $T^{(t)}_z$ by the Markov property and the fact that $0$ is a
regular point for L\'evy processes with infinite total variation.
Therefore, writing $|I^t_i|$ for the Lebesgue measure of $I^t_i$,  
$T_z-\sum_k' |I^t_k|$ (where the sum is over the $I^t_k$ that are maximal) has same law as $T^{(t)}_z$, and in particular it is
nonzero a.s. Now to check that the intervals
$I^t_i$ are separated (we are only interested by those which are maximal), 
consider two left-ends of such intervals such as
$\tau^{t,a}_i<\tau^{t,a}_j$ (where $a$ is small enough). The regularity of
$0$ for the L\'evy process $X$ implies that
$\inf_{s\in[\sigma^{t,a}_i,\tau^{t,a}_i]}X_s< X_{\sigma^{t,a}_i}$, so
by the same arguments as above and the Markov property at
$\sigma^{t,a}_i$, there exists a (random)
$\eps^a_{i,j}>0$ such that given $\eps^a_{i,j}$,
$$\tau^{t,a}_j-\sigma^{t,a}_i-\sum_{I^t_k\subset
[\sigma^{t,a}_i,\tau^{t,a}_j]}|I^t_k|\ind_{\{I^t_k \mbox{ maximal}\}}$$
is stochastically larger than $T^{(t)}_{\eps^a_{i,j}}$. 
This ensures the a.s.\ separation of the
$I^t_k$'s, so the a.s.\ convergence of
${\tt UNPLUG}(X,(I^t_i:z^t_i>a))$ to ${\tt UNPLUG}(X,(I^t_i,i\geq 1))$ as
$a\downarrow0$ comes from Lemma \ref{unpluglmm}. Identifying the limiting
law follows from the above discussion.
\cq

Now let as before
$X^t_i(s)=X_{\tau^t_i+s}-X_{\tau^t_i}$ for $0\leq s\leq \sigma^t_i-
\tau^t_i$ and
$i\geq 0$, where by convention $\tau^t_0=0$, and $\sigma^t_0=T_1$. We write
$-\tau^t_i+I^t_k=[\tau^t_k-\tau^t_i,\sigma^t_k-\tau^t_i]$. The next lemma does most of
the job to extract the different tree components of the logged stable
tree at time $t$.  

\begin{lmm}\label{unpluglem}
{\rm (i)} Under the law $\bP^{(-1,\infty)}$, as $a\downarrow0$, the
processes 
${\tt UNPLUG}(X^t_i,(-\tau^t_i+I^t_k,k:I^t_k\subsetneq I^t_i
\mbox{ and }z^t_k>a)),i\geq 1$ converge in $\D^{\dagger}([0,\infty),\R)$  
to the processes 
$Y^t_i={\tt UNPLUG}(X^t_i,(-\tau^t_i+
I^t_k,k:I^t_k\subsetneq I^t_i)),i\geq1$. 

{\rm (ii)} The process $Y^t_i$ has the same law as $z^t_i+X^{(t)}$ under 
$\bP$, killed when it first hits $0$, 
and these processes are independent conditionally on $(z^t_i,i\geq 1)$.

{\rm (iii)} The sum of the durations of $Y^t_i,i\geq 0$ equals $T_1$ a.s.
\end{lmm}


\proof
(i) Fix $a>0$, we modify slightly the notations of the preceding proof by 
letting $\tau^{t,a}_1<\ldots<\tau^{t,a}_{k(a)}$ be the times when 
$Z^{(t)}$ accomplishes a jumps that is $>a$, and letting $\sigma^{t,a}_i=
\inf\{u\geq \tau^{t,a}_i:X_u=X_{\tau^{t,a}_i-}\}$. 
Let also $\tau^{t,a}_0=0, \sigma^{t,a}_0=T_1$. 
Write $I^{t,a}_i=[\tau^{t,a}_i,\sigma^{t,a}_i]$, and
let $X^{t,a}_i(s)=X_{\tau^{t,a}_i+s}-X_{\tau^{t,a}_i}$ for 
$0\leq s\leq \sigma^{t,a}_i-\tau^{t,a}_i$. 
By the Markov property at times $\tau^{t,a}_i,\sigma^{t,a}_i$, 
we obtain that for every $i$, $X^{t,a}_i$ is independent
of ${\tt UNPLUG}(X,I^{t,a}_i)$ given the 
jump $\Delta X_{\tau^{t,a}_i}$. By a repeated use of the Markov property, 
we obtain the independence of the processes ${\tt UNPLUG}(X^{t,a}_i,
(-\tau^t_i+I^{t,a}_k:I^{t,a}_k\subsetneq I^{t,a}_i))$ given 
$(\Delta X_{\tau^{t,a}_i},1\leq i\leq k(a))$, and moreover, the law of 
$X^{t,a}_i$ given 
$\Delta X_{\tau^{t,a}_i}$ is that of $X$ under $P$, killed when it first 
hits $-\Delta X_{\tau^{t,a}_i}$.
Letting $a\downarrow 0$ and applying Lemma \ref{plugg} finally 
gives the convergence to the processes $Y^t_i$, as well as
the conditional independence and the distribution
of the processes, giving also (ii). 

(iii) 
Let us introduce some extra notation. Say that the marked 
jump with magnitude $z^t_i$ is of the $j$-th kind if and only if
the future infimum process 
$(\inf_{s\leq u\leq \tau^t_i} X_u,0\leq s\leq \tau^t_i)$ accomplishes
exactly $j$ jumps at times that correspond to marked jumps of $X$. 
Write $|I^t_i|$ for the duration of $X^t_i$ and let $A_j$ be 
the set of indices $i$ such that $\tau^t_i$ is a jump time of the 
$j$-th kind. By a variation of Lemma \ref{finitehub} already used above, 
every marked jump is of the $j$-th kind 
for some $j$ a.s. By Lemma \ref{plugg} the duration of $Y^t_0$ is 
$T_1-\sum_{i\in A_1}|I^t_i|$, similarly, one has that if $i\in
A_j$,  the duration of $Y^t_i$ equals
$|I^t_i|-\sum_{k\in A_{j+1}}|I^t_k|\ind_{\{I^t_k\subset I^t_i\}}$. 
Therefore, proving that the 
sum of durations of $Y^t_i$ equals $T_1$ amounts to showing that 
$\sum_{i\in A_j}|I^t_i|\to 0$ in probability as $j\to \infty$. 
But the sum of the marked jumps is finite a.s., since conditionally 
on a marked jump $z^t_i$, the duration of the corresponding $X^t_i$
has same law as $T_{z^t_i}$, and since we have independence as $i$
varies. Hence this sum is (conditionally on $(z^t_i,i\geq 1)$) 
equal in law to $T_{\sum_{i\in A_j} z^t_i}$ under $\bP$, and it converges 
to $0$. \cq

\begin{lmm}\label{sfragplus}
The process $(F^+(t),t\geq 0)$ is a Markovian self-similar fragmentation 
with index $1/\alpha$. Its erosion coefficient is $0$
\end{lmm}

\proof
For every $v>0$, define the processes $X^t_i$ under $\bN^{(v)}$ as in 
the preceding section, replacing the duration $1$ by $v$. 
By virtue of Lemma \ref{unpluglem} and by excursion theory, 
we obtain that for almost every $v>0$, and for all $t$ in a dense
countable  subset of $\R_+$, under $\bN^{(v)}$, the processes
${\tt UNPLUG}(X^t_i,(-\tau^t_i+I^t_k:I^t_k\subsetneq I^t_i\mbox{ and }z^t_k>a))$
converge as $a\downarrow0$ to processes $Y^t_i$ that are independent
conditionally  on the $z^t_i$'s and on their durations, and whose
durations sum to $v$ (by convention we let $X^t_0=X$). By scaling, this
statement remains  valid for $v=1$. We then extend it to all $t\geq 0$ by
a continuity argument. The case $t=0$ is obvious, so take
$t_0>0$ and $t\uparrow t_0$ in the dense subset of $\R_+$. Almost
surely, $t_0$ is not a time at which a new hub is marked, so
$X^{t_0}_i=X^t_i$ for $t$ close enough of $t_0$, and by Lemma
\ref{unpluglmm} and the fact that
$\{I^t_i,i\geq 0\}\subset\{I^{t_0}_i,i\geq 0\}$ for $t\leq t_0$,
$$Y^{t_0}_i={\tt UNPLUG}(X^t_i,(-\tau^t_i+I^{t_0}_k:I^{t_0}_k\subsetneq
I^{t_0}_i))=\lim_{t\uparrow t_0}{\tt UNPLUG}(X^t_i,(-\tau^t_i+I^t_k:I^t_k\subsetneq
I^t_i)).$$
Now recall the notation $X^t_{i,j},I^t_{i,j}=[\tau^t_{i,j},\sigma^t_{i,j}]$ 
from Sect.\ \ref{fplusalt}, and for $j\geq 1$ write 
$Y^t_{i,j}={\tt UNPLUG}(X^t_{i,j},(-\tau^t_{i,j}+I^t_k:I^t_k\subsetneq
I^t_{i,j}))$ 
for the excursions of $Y^t_i$ above its infimum, 
ranked in the order corresponding to $X_{i,j}$. Then by the same arguments as 
in the proof of Lemma \ref{plugg}, the joint law of the
durations of 
$Y^t_0,Y^t_{i,j},i\geq 1,j\geq 1$ equals the law of 
$(|C^t_0|,|C^t_{i,j}|,i\geq 1,j\geq 1)$ with 
notations above. 
Hence, by Lemma \ref{cij} and the fact that 
excursions of $X^{(t)}$ with prescribed duration are stable excursions, it 
holds that
conditionally on $F^+(t)=(x_1,x_2,\ldots)$, the excursions $Y^t_{i,j}$ are
independent stable excursions with respective durations
$x_1,x_2,\ldots$. 
 
Now let $\sim^{t,i,j}_{t'}$ be the equivalence relation
defined for the excursion $Y^t_{i,j}$ in a similar way as $\sim_t$ for the
normalized  excursion of $X$. Write also
$j_t(u)=u-\tau^t_{i,j}-\sum_{k:I^t_k\subsetneq I^t_i,\sigma^t_k<u}|I^t_k|$ 
for $u\in[0,1]$, whenever $u\in C^t_{i,j}$. 
Then it is clear that if $x,y\in C^t_{i,j}$,
one has also $x\sim_{t+t'}y$ if and only if
$j_t(x)\sim^{t,i,j}_{t'}j_t(y)$. By the scaling property, a
stable  excursion $\eps^x$ with 
duration $x$ where every jump with magnitude $\ell$ is marked with 
probability $1-\exp(-t' \ell)$ is obtained by taking a normalized excursion 
$(\eps^{1}_s,0\leq s\leq 1)$, marking every jump with magnitude $\ell$
independently with probability $1-\exp(-t' x^{1/\alpha}\ell)$, and then 
letting $\eps^x_s=x^{1/\alpha}\eps^1_{s/x}$ for $0\leq s\leq x$; the 
marked jumps of $\eps^x$ occurring at the times $sx$ whenever $s$ is 
a marked jump time for $\eps^1$. This means that given $F^+(t)=(x_1,\ldots)$, 
the process $(F^+(t+t'),t'\geq 0)$ has the same law as 
$((x_1F^{+,1}(x_1^{1/\alpha}t'),x_2 F^{+,2}(x_2^{1/\alpha}t'),
\ldots)^{\downarrow},
t'\geq 0)$ where the $F^{+,i}$'s are independent copies of 
$F^+$. This entails both the Markov property and the self-similar property, 
the self-similarity index being $1/\alpha$. 
Moreover, Lemma \ref{unpluglem} (iii) shows that the sum of durations of
$Y^t_{i,j}$ is $1$ a.s.\ under
$\bN^{(1)}$, so $\sum_i F^+_i(t)=1$ a.s.\ and 
the erosion coefficient must be $0$ according to \cite{bertsfrag02}. 

To conclude, we notice that the previous result of continuity
in probability of $F^+$ at time $0$ extends to any time $t\geq 0$ by
the self-similar fragmentation property. 
\cq

\subsection{Splitting rates and dislocation measure}

To complete the study of the characteristics of $F^+$, we must 
identify the dislocation measure. This is done by computing 
the {\em splitting rate} of 
the stable tree, that is, the rate at which the tree with mass $1$ 
instantaneously splits into a sequence of subtrees with given masses
$s_1\geq s_2,\ldots$ with $\sum_i s_i=1$, 
by analogy with the splitting rate of the Brownian CRT in 
\cite{jpda98sac}. 

We will need the following lemma from \cite{MS03}, which is similar to
L\'evy's method to compute the jump measure of a L\'evy process. 

\begin{lmm}\label{msprp}
Let $(F(t),t\geq 0)$ be a self-similar fragmentation with index
$\beta\geq 0$ and erosion coefficient $c=0$. Then for every function $G$
that is continuous and null on a neighborhood of $(1,0,\ldots)$ in
$S$, 
$$t^{-1}E[G(F(t))]\build\to_{t\downarrow0}^{}\nu(G).$$
\end{lmm}

Recall that our marking process on the 
hubs of the tree amounts to taking a Poisson process with intensity 
$m(\d v)=\sum_b L(b)\delta_b(\d v)$ 
on $\TT$, where the sum is over hubs $b\in \TT$. For $v\in \TT$, let
$\TT_1(v),\TT_2(v),\ldots$ be the tree components of the forest 
obtained when removing $v$, arranged by decreasing order of masses, and let 
$$r(\d\bs)=N^{(1)}\left(
m\{v\in \TT:(\mu(\TT_1(v)), \mu(\TT_2(v)),\ldots)\in \d\bs\}\right)$$
be the rate at which a $m$-picked vertex splits $\TT$ into trees with 
masses in  a volume element  $\d\bs$ (recall that  the stable tree  is defined
under the normalized excursion law $N^{(1)}$). It is quite intuitive that 
the splitting rate equals the dislocation measure of $F^+$, 
and Theorem \ref{T1} reduces to the two 
following lemmas:

\begin{lmm}\label{sprate}
The splitting rate $r(\d\bs)$ equals the dislocation measure $\nu_+$ of $F^+$. 
\end{lmm}

\proof
For $t\geq0$ we let $\TT(t)$ be the forest obtained by our logging
procedure of the stable tree at time $t$. Let $n\geq 2$, and consider $n$
leaves $L_1,\ldots,L_n \in \TT$ that are independent and distributed
according to the mass measure
$\mu$, conditionally on $\mu$ (we are implicitly working on an enlarged 
probability space). Write $\Pi_n(t)$ for the partition of
$[n]=\{1,\ldots,n\}$ obtained by letting $i$ and $j$ be in the same
block of $\Pi_n(t)$ if and only if $L_i$ and $L_j$ belong to the same
tree component of $\TT(t)$. For $K>2$ let $\Lambda^n_K(t)$ 
be the event that at time $t$, the leaves $L_1,\ldots,L_n$
are all contained in tree components of $\TT(t)$ with masses $>1/K$. 
Write ${\cal P}^*_n$ for the set of partitions $\pi$ of 
$[n]=\{1,\ldots,n\}$ with at least two non void blocks
$A_1,\ldots,A_k$ (for some arbitrary ordering
convention). 
Given $F^+(t)=\bs=(s_1,s_2,\ldots)$, the probability that $\Pi_n(t)$ equals 
some partition $\pi\in{\cal P}^*_n$ and that $\Lambda^n_K(t)$ happens is 
$$G_K(\bs)=\bN^{(1)}(\Pi_n(t)=\pi, \Lambda^n_K(t)|F^+(t)=\bs)=
\sum^{*K}_{i_1,\ldots,i_k}\prod_{j=1}^k s_{i_j}^{\#A_j},$$
the sum being over pairwise distinct $i_j$'s such that $s_{i_j}>1/K$.
This last function is continuous and
null on a neighborhood of $(1,0,\ldots)$, so 
Lemma \ref{msprp} (which we may use by Lemma \ref{sfragplus})
gives 
\begin{equation}\label{expnu}
\lim_{t\downarrow 0} t^{-1}\bN^{(1)}(\Pi_n(t)=\pi,\Lambda^n_K(t))=\int_S
\nu_+(\d\bs)\sum_{i_1,\ldots,i_k}^{*K}\prod_{j=1}^k s_{i_j}^{\#A_j}.
\end{equation}
We claim that knowing this quantity for every $n,\pi,K$ characterizes $\nu_+$.
One can obtain this by first letting 
$K\to \infty$ by monotone convergence, and then using an argument 
based on exchangeable partitions as in \cite[p.\ 378]{kingman78} 
(a Stone-Weierstrass argument can also work). 

On the other hand, for any $b$ in the set ${\cal H}(\TT)$ of branchpoints of 
$\TT$, let $\pi^b_n$ be the 
partition of $[n]$ obtained by letting $i$ and $j$ be in the same block if 
and only if $b$ is not on the path from $L_i$ to $L_j$. Let 
also $\TT_{L_i}(b)$ be the tree component of the forest obtained by 
removing $b$ from $\TT$ that contains $L_i$. 
For $K\in(2,\infty]$ and $\pi\in{\cal P}^*_n$, let
$\Psi^n_K(\pi)$ be the set of branchpoints $b\in\TT$ such that 
$\pi^b_n=\pi$ and such that $\mu(\TT_{L_i}(b))>1/K$ for $1\leq i\leq n$,
and let $\Psi^n_K=\bigcup_{\pi\in{\cal P}^*_n}\Psi^n_K(\pi)$. Recall that 
we may construct the fragmentation $F^+$ by cutting the stable tree
at the points of a Poisson point process 
$(b(s),s\geq 0)$ with intensity $\d s\otimes m(\d b)$. 
Now for $\Pi_n(t)=\pi$ to happen, 
it is plainly necessary that at least one $b(s)$ falls in $\Psi^n_{\infty}$ 
for some $s\in[0,t]$, if in addition $\Lambda^n_K(t)$ happens then no 
$b(s),0\leq s\leq t$ must fall in $\Psi^n_{\infty}\setminus \Psi^n_K$. 
Therefore, 
\begin{equation}\label{intermexp}
\bN^{(1)}(\Pi_n(t)=\pi,\Lambda^n_K(t))=\bN^{(1)}\left(\exists!\, 
s\in[0,t]:b(s)\in \Psi^n_{\infty}, 
\mbox{ and }b(s)\in \Psi^n_K(\pi),\Lambda^n_K(t)\right)+R(t),
\end{equation}
where the residual $R(t)$ is bounded by the probability that
$b(s)$ falls in $\Psi^n_{\infty}$ for at least two $s\in[0,t]$.
Hence $R(t)=o(t)$ by standard properties of Poisson processes 
provided we can show that $N^{(1)}[m(\Psi^n_{\infty})]<\infty$. This could be shown 
using the forthcoming lemma, but we may also just notice that if 
$N^{(1)}[m(\Psi^n_{\infty})]$ was infinite, then there would be arbitrarily 
many $b(s),0\leq s\leq t$ 
falling in $\Psi^n_{\infty}\setminus\Psi^n_K$ for some appropriately
large $K$, and the probability in (\ref{intermexp}) would be $0$, which is impossible from 
the beginning of this proof and since $F^+$ is a self-similar
fragmentation with nonzero dislocation measure (because it has erosion
coefficient $0$ and it is not constant). On the other hand,  conditionally
on the event on the right-hand side of (\ref{intermexp}),  the $b(s),0\leq
s\leq t$ that do not fall in
$\Psi^n_{\infty}$ (call them 
$b'(s)$) form an independent Poisson point 
process with intensity $m(\cdot\cap{\cal H}(\TT)\setminus \Psi^n_{\infty})$. 
Therefore, the size of the tree component of the forest obtained when 
removing the points $b'(s),0\leq s\leq t$ that contains $L_1$ converges 
a.s.\ to $1$ as $t\downarrow0$ (so it also contains the other $L_i$'s for 
small $t$ a.s.), as it is stochastically bigger than the 
component of $\TT(t)$ containing $L_1$, and since
$F^+(t)\to(1,0,\ldots)$ in probability as
$t\downarrow 0$. It follows that one can remove $\Lambda^n_K(t)$
from the right-hand side of (\ref{intermexp}), and basic properties 
of Poisson measures finally give 
$t^{-1}\bN^{(1)}(\Pi_n(t)=\pi,\Lambda^n_K(t))\to 
\bN^{(1)}[m(\Psi^n_K(\pi))]=N^{(1)}[m(\Psi^n_K(\pi))]$. This  last quantity is
finally equal
to $\int_S r(\d\bs)\sum_{i_1,\ldots,i_k}^{*K}\prod_{j=1}^k s_{i_j}^{\#A_j}$
since $L_i$ belongs to $B\subset\TT$ with probability $\mu(B)$ that is 
equal to the Lebesgue measure of the subset of $[0,1]$ encoding $B$. 
Identifying with (\ref{expnu}) gives the claim. \cq

\begin{lmm}\label{rform}
One has $r(\d\bs)=\nu_{\alpha}(\d\bs)$ with the notations of Theorem \ref{T1}.
\end{lmm}

\proof
We must see what is the effect of splitting $\TT$ at a hub $b$ picked 
according to $m(\d v)$. By definition, $m$ picks a hub proportionally to 
its local time, and by Proposition \ref{onetoone}, 
hubs are in one-to-one correspondence with jumps
of the stable excursion with duration $1$. More precisely, if 
$b$ is the hub that has been picked and with the notations 
$\tau(b),\sigma(b)$ above, the masses of the tree components obtained 
when removing $b$ are equal to the lengths of the constancy intervals
of the infimum process of $(X_{\tau(b)+s}-X_{\tau(b)},0\leq s\leq \sigma(b)
-\tau(b))$, and the extra term $1-(\sigma(b)-\tau(b))$. By Vervaat's 
theorem, we may suppose that the excursion is the Vervaat transform of a
stable bridge and that the marked jump in the excursion
corresponds to a jump $(s,\Delta X_s)$ of the bridge picked according
to the $\sigma$-finite measure $\sum_{u:\Delta X_u>0}\Delta X_u 
\delta_{(u,\Delta X_u)}(\d s,\d x)$.  By Lemma \ref{markbr}, this marked jump 
equals $(s,x)$ according to a certain $\sigma$-finite ``law'', while given 
$(s,\Delta X_s)=(s,x)$, the bridge $X$ 
has the same law as $X\oplus(s,x)$, under the law 
$P^1_{0\to-x}$.  

Therefore, we have have obtained a representation of the
excursion together with a marked jump as a bridge
$X$ with law $P^1_{0\to -x}$, where $x$ is independent with
some $\sigma$-finite ``law'', to which has been added the marked jump of size 
$x$ at an independent uniform time $s$, and which has finally undergone the
Vervaat transformation. Using the invariance of bridge laws under independent 
cyclic shifts, it is now easy to see that the lengths of the constancy 
intervals of $(X_{\tau(b)+s}-X_{\tau(b)},0\leq s\leq\sigma(b)-\tau(b))$ 
defined above have the same law as the intervals of 
constancy of the infimum process of $(X_{s'+s}-X_s,0\leq s'\leq T_x)$ under 
$P^1_{0\to-x}$ (with $x$ as above), while the remaining term 
$1-(\sigma(b)-\tau(b))$ has (jointly) law $1-T_x$. 

It is now easy that conditionally on 
$x,T_x=t$ these constancy intervals have the same law as 
$\Delta T_{[0,x]}$ given $T_x=t$ under $P$ (one actually checks that 
$(X_u,0\leq u\leq T_x)$ is the first-passage bridge with law 
$P^t_{0\downarrow -x}$ defined before Lemma \ref{tatb} below). 
The law of $1-T_x$ given $x$ is simply obtained by using the definition 
of bridges and the Markov property: for $a<1$ and positive measurable $f$, 
\begin{eqnarray*}
E^1_{0\to -x}[f(1-T_x)\ind_{\{T_x<a\}}]&=&E^1[f(1-T_x)\ind_{\{T_x<a\}}
p_1(-x)^{-1}p_{1-a}(-x-X_a)]\\
&=&\int_0^a \d s\, q_x(s)f(1-s) 
p_1(-x)^{-1}\int \d y\, p_{a-s}(y) p_{1-a}(-y)\\
&\build\to_{a\to 1}^{}&\int_0^1 \d s\,  q_x(s)f(1-s)p_{1-s}(0)p_1(-x)^{-1}.
\end{eqnarray*}
In the last integral, change variables $1-s\to s$, use $p_1(-x)=x^{-1}q_x(1)$,
check by scaling that $p_s(0)=s^{-1/\alpha}p_1(0)$, and 
conclude by identifying with Lemma \ref{permanform} that 
$1-T_x$ under $P^1_{0\to-x}$ has same law as a size-biased pick from 
$\Delta T_{[0,x]}$ given $T_x=1$ under $P$ (notice that in particular 
we must have $p_1(0)=c_{\alpha}$). By Lemma \ref{permanform} (ii), it follows
that given the local time $x$ of the marked hub $b$, 
the law of the sizes of the stable tree split at this hub is the same as that 
of $\Delta T_{[0,x]}$ given $T_x=1$ under $P$. 

Putting pieces together and recalling the distribution of the marked 
jump $x$ from Lemma \ref{markbr} we obtain the formula 
$$r(\d\bs)=\int_0^{\infty}\d x \frac{C_{\alpha}p_1(-x)}{x^{\alpha}p_1(0)}
P(\Delta T_{[0,x]}\in \d\bs|T_x=1).$$
By using the scaling property for $T$ and its density ($q_x(1)=x^{-\alpha}q_1
(x^{-\alpha})$), formula (\ref{ballotc}) and a change of 
variables, we obtain 
\begin{eqnarray*}
r(\d\bs)&=& \int_0^{\infty}\d x\frac{C_{\alpha}q_1(x^{-\alpha})}{
c_{\alpha}x^{2\alpha+1}}
P(T_1^{-1}\Delta T_{[0,1]}\in \d\bs|T_1=x^{-\alpha})\\
&=&\alpha^{-1}c_{\alpha}^{-1}C_{\alpha}\int_0^{\infty}\d u\, u\,
q_1(u)P(T_1^{-1}
\Delta T_{[0,1]}\in \d\bs|T_1=u),
\end{eqnarray*}
which gives the desired formula, after checking that 
$\alpha^{-1}c_{\alpha}^{-1}C_{\alpha}=D_{\alpha}$. \cq

\section{Study of $F^{\natural}$}\label{levyrep}

Recall the construction of $F^{\natural}$ (under the measure $\bN^{(1)}$)
from Sect.\ \ref{intro}. As noticed above, this fragmentation process 
somehow generalizes the one
considered in \cite{bertfrag99,mier01} (we could actually build it in an
analogous way for a large class of L\'evy processes with no negative
jumps, though the resulting fragmentations would not be self-similar due
to the absence of scaling). Notice that none of the fragmentation processes
of \cite{mier01} are self-similar, but for the Brownian case. The reason for
this was a lack of a Girsanov-type theorem saying that a L\'evy process
plus drift has a law that is absolutely continuous with the initial
process, but for the Brownian case. Here, this is fixed by Proposition
\ref{abscont}, but where the operation is removing jumps rather than
adding a drift. 

\subsection{The self-similar fragmentation property}

For any $t'>t\geq 0$ let $\mu_t(x,\d\bs)$  be a kernel from $\R_+^*$ to $S$
defined as follows: $\mu_t(x,\d\bs)$ 
is the law of the ranked lengths of the constancy intervals of the
process $\un{X}^{(t)}$ under $\bN^{(x)}$. Moreover, define
$F^{\natural,1}$ exactly as $F^{\natural}$, but where $X$ is under the law
$\bP^{(-1,\infty)}$. In particular, $F^{\natural,1}(t)$ is not $S$-valued
(the sum of its components is random). 

\begin{prp}\label{elem}
{\rm (i)} The processes $F^{\natural,1}$ and $F^{\natural}$ 
enjoy the fragmentation property, with 
fragmentation kernel
$\mu_t(x,\d\bs)$. That is, conditionally on $F^{\natural,1}(t)
=(x_1,x_2,\ldots)$
(resp.\ $F^{\natural}(t)$), $F^{\natural,1}(t+t')$ 
(resp.\ $F^{\natural}(t+t')$) has the same law as the decreasing
rearrangement of independent sequences $\bs_i$ with respective laws
$\mu_{t'}(x_i,\d\bs)$. 

{\rm (ii)} The process $F^{\natural}$ is a self-similar fragmentation 
with index $1/\alpha$, and no erosion. 
\end{prp}

The fact that $F^{\natural}$ 
is a fragmentation process directly comes from the fact that
the processes $X^{(t')}-X^{(t)}=Z^{(t)}-Z^{(t')}$ are non-increasing. We now
prove the fragmentation property. The  key lies in a Skorokhod-like relation 
that is analogous to that in \cite{bertfrag99} and generalized in
\cite{mier01}.

\begin{lmm}\label{skoform}
For every $t,t'\geq 0$ and $s\geq 0$, one has
$$\un{X}^{(t)}_s=\inf_{0\leq u\leq
s}(\un{X}^{(t+t')}_u+(Z^{(t+t')}_u-Z^{(t)}_u)).$$
\end{lmm}

The proof can be done following exactly the same lines as in
\cite[Lemma 2]{bertfrag99}. As a consequence, we obtain that the sigma-field
${\cal G}_t=\sigma\{X^{(t)},(Z^{(s)},0\leq s\leq t)\}$ induces a filtration,
with respect to which $F^{\natural,1}$ is adapted.  

The end of the proof of the fragmentation property in Proposition \ref{elem}
also goes as in \cite{bertfrag99}. For any variable $K$ that is ${\cal
G}_t$-measurable, the excursions of $X^{(t)}$ above its infimum and before time
$T^{(t)}_K$ are independent excursions conditionally on ${\cal G}_t$,
respectively conditioned to have durations
$\ell^{(t)}_{1,X},\ell^{(t)}_{2,X},\ldots$ where the last family is the
decreasing sequence of constancy intervals of $\un{X}^{(t)}$ before time
$T^{(t)}_{K}$. Take $K=\un{X}^{(t)}_{T_1}$, which is measurable with respect
to ${\cal G}_t$ by virtue of the Skorokhod property. Then $T^{(t)}_K=T_1$,
which gives readily that conditionally on ${\cal G}_t$, the excursions of
$X^{(t)}$ above $\un{X}^{(t)}$ are independent with durations 
$(F^{\natural,1}_i(t),i\geq 0)$.

To conclude, it remains to notice that the lack of memory of the exponential
law implies that the jumps that are unmarked at time $t$ but that are marked
at time $t+t'$ can be obtained also by marking with probability
$1-e^{-t'\ell}$ any unmarked jump at time $t$ that has magnitude
$\ell$. Thus, conditionally on $F^{\natural,1}(t)$, 
we obtain a sequence with the same
law as $F^{\natural,1}(t+t')$ 
by taking independent sequences $(\bs_i,i\geq 1)$ with laws
$\mu_{t'}(F^{\natural,1}_i(t),\d\bs)$ and rearranging, as claimed. 
This remains true for $F^{\natural}$ by excursion
theory and scaling. 

To show the self-similarity for $F^{\natural}$, it then suffices to check, 
using the scaling property of the excursions of stable processes, that 
$\mu_t(x,\d\bs)$ is the image of $\mu_{tx^{1/\alpha}}(1,\d\bs)$ by 
$\bs\mapsto x\bs$. The fact that $F^{\natural}$ has no erosion again comes
from the fact that $\sum_i F^{\natural}_i(t)=1$ a.s.

\subsection{The semigroup}

According to the preceding section, and since plainly there is no loss 
of mass in the fragmentation $F^{\natural}$ (so the erosion coefficient is 
$0$), proving Theorem \ref{T2} requires only 
to check that the dislocation measure of $F^{\natural}$ 
equals that of $F^+$. It is intuitively straightforward that this is the 
case, by looking at the procedure we use for deleting jumps, and indeed we 
could easily follow the same lines as above and compute a ``splitting rate''
for the bridge, when the ``first'' marked jump is deleted. However, a
nice feature of this fragmentation is that we can compute explicitly its
semigroup (hence that of $F^+$), as will follow. The semigroup then gives
enough information to re-obtain the dislocation measure, and this will
prove Theorem \ref{T2}. Recall from Sect.\ \ref{levyproc} that
$\rho^{(t)}_1$ is the density of $Z^{(t)}_1$ under $\bP$.

\begin{prp}\label{prpsemig}
The semigroup of $F^{\natural}$ is given by 
$$\bN^{(1)}(F^{\natural}(t)\in \d\bs)=
\int_0^{\infty}\d z\frac{p_1^{(t)}(-z)\rho_1^{(t)}(z)}{p_1(0)}
P(\Delta T_{[0,z]}\in \d\bs| T_z=1).$$
\end{prp}

We will need a couple of intermediate lemmas. 
Since $Z^{(t)}$ is non-decreasing, under the law $\bN^{(1)}$, the process 
$X^{(t)}$ starts at $0$ and hits $-Z^{(t)}_1$ at time $1$ for the first time. 
Since we are interested in the constancy intervals of $\un{X}^{(t)}$, 
and thanks
to Vervaat's theorem, we would like to relate these constancy intervals 
to the bridge of $X$. We now work under the law of the bridge 
with unit duration $\bP_{0\to0}^1$, so we may suppose that the 
excursion of $X$ with duration $1$ is equal to the Vervaat transform
$VX$. Let $m=-\un{X}_1$ be the absolute value of the 
minimum of $X$, and $\tau_2=T_{m-}$ be the 
(a.s.\ unique) time when $X$ attains this minimum, so $VX=V^{\tau_2}X$. 
Decompose 
$X$ as $X^{(t)}+Z^{(t)}$ where $Z^{(t)}$ is the cumulative process of marked 
jumps. Then $VX=V^{\tau_2}X^{(t)}+V^{\tau_2}Z^{(t)}$, and by independence of 
the marking procedure of jumps we can consider that $V^{\tau_2}Z^{(t)}$ is the 
cumulative process of marked jumps for the excursion $VX$. The problem is 
now to describe the law of lengths of the constancy intervals of the 
process $\un{V^{\tau_2}X}^{(t)}$. Let $m^{(t)}=
-\un{X}^{(t)}_1$ be the absolute value of the minimum of $X^{(t)}$
and $\tau_3=T^{(t)}_{m^{(t)}-}$ be the (a.s.\ unique) time when 
$X^{(t)}$ attains this minimum. Let also $\tau_1=T^{(t)}_{m^{(t)}-Z^{(t)}_1}$
be the first time when $X^{(t)}$ attains the value $Z^{(t)}_1-m^{(t)}$. 
The following lemma is somehow ``deterministic''. 
For $a<b$, write $X_{[a,b]}$ for the process $(X_{a+s}-X_a,0\leq s\leq b-a)$.


\begin{lmm}
One has $\tau_1\leq \tau_2\leq \tau_3$ a.s., and 
the sequence of 
lengths of the constancy intervals of $\un{V^{\tau_2}X}^{(t)}$, 
ranked in decreasing order, is equal to that of the process 
$\un{X}^{(t)}_{[\tau_1,\tau_3]}$, to which has been added (at the appropriate
rank) the extra term $1-\tau_3+\tau_1$.
\end{lmm}

\proof
Since $Z^{(t)}$ is an increasing process, one has 
$X^{(t)}_{\tau_2}=X_{\tau_2}-Z^{(t)}_{\tau_2}\leq X_s-Z^{(t)}_s$ for any 
$s\leq \tau_2$. Hence, $X^{(t)}_{\tau_2}=\un{X}^{(t)}_{\tau_2}$ which 
implies $\tau_2\leq \tau_3$. On the other hand, one has 
$-m^{(t)}=X_{\tau_3}-Z^{(t)}_{\tau_3}\geq -m-Z^{(t)}_1$ and 
thus $m^{(t)}-Z^{(t)}_1\leq m$, implying $\tau_1\leq \tau_2$. 

For convenience, if $(f(x),0\leq x\leq \zeta)$ and $(f'(x),0\leq x\leq \zeta')$
are two c\`adl\`ag functions, we let $f\bowtie f'$
be the concatenation of the paths of $f$ and $f'$, defined by 
$$f\bowtie f'(s)=\left\{\begin{array}{cl}
f(s) & \mbox{if } 0\leq s< \zeta\\
f'(s-\zeta)+f(\zeta) &\mbox{if }\zeta\leq s\leq \zeta+\zeta'
\end{array}\right. .$$ 
We let 
$Y^1=X^{(t)}_{[0,\tau_2]}$, $Y^2=X^{(t)}_{[\tau_2,\tau_3]}$ and 
$Y^3=X^{(t)}_{[\tau_3,1]}$, so $X^{(t)}=Y^1\bowtie Y^2\bowtie Y^3$, 
and $V^{\tau_2}X^{(t)}=Y^2\bowtie Y^3\bowtie Y^1$.

Observing that $Y_3$ is non-negative, we obtain that 
$\un{Y^2\bowtie Y^3}=\un{Y}^2\bowtie {\bf 0}_{[0, 1-\tau_3]}$ where 
${\bf 0}_{[0,a]}$ is the null process on $[0,a]$. Since the final value 
of $Y_3$ is $m^{(t)}-Z^{(t)}_1$, we obtain that 
$$\un{V^{\tau_2}X}^{(t)}=\un{Y}^2\bowtie {\bf 0}_{[0, 1-\tau_3]}\bowtie
{\bf 0}_{[0,\tau_1]}\bowtie \un{X}^{(t)}_{[\tau_1,\tau_2]}
=\un{Y}^2\bowtie{\bf 0}_{[0,1-\tau_3+\tau_1]}\bowtie \un{X}^{(t)}_{
[\tau_1,\tau_2]}.$$
It follows that the constancy intervals of $\un{V^{\tau_2}X}^{(t)}$ are 
the same as those of $\un{X}^{(t)}$, except for the first and last constancy 
intervals of $\un{X}^{(t)}$ which are merged to form the 
constancy interval with length $1-\tau_3+\tau_1$. \cq

The rest of the section is devoted to the study of these constancy intervals. 
Recall from Lemma \ref{decompX} 
that under $\bP^1_{0\to0}$, the process $X^{(t)}$ has law
$P^1_{0\to -{\cal Z}}$, where ${\cal Z}$ is an independent random variable 
with law $P({\cal Z}\in \d z)=p_1(0)^{-1}p^{(t)}_1(-z)\rho^{(t)}_1(z)\d z$. 
It thus suffices to analyze the constancy intervals of 
$\un{X}_{[\tau_1,\tau_3]}$ under 
the law $P^1_{0\to -z}$ for fixed $z>0$, where we now call $m=-\un{X}_1$, 
$\tau_1$ the time when $X$ first hits level $z-m$ and 
$\tau_3$ the first time when $X$ attains level $-m$. 

For $z>0$, let $(P^v_{0\downarrow -z},v>0)$ be 
a regular version of the conditional law $P^{(-z,\infty)}[\cdot| T_z=v]$. 
Call this the law of the first-passage bridge from $0$ to $-z$ with 
length $v$. A consequence of the Markov property is

\begin{lmm}\label{tatb}
Let $a,b>0$. For (Lebesgue) almost every $v>0$, 
under the law $P^v_{0\downarrow-(a+b)}$, the law 
of $T_a$ is given by 
$$P^v_{0\downarrow-(a+b)}(T_a\in \d s)=\d s\frac{q_a(s)q_b(v-s)}{
q_{a+b}(v)}.$$
Moreover, conditionally on $T_a$, the paths $(X_s,0\leq s\leq T_a)$ and 
$(X_{s+T_a}-a,0\leq s\leq T_{a+b}-T_a)$ are independent with respective laws 
$P^{T_a}_{0\downarrow -a}$ and $P^{v-T_a}_{0\downarrow -b}$.
\end{lmm}

We also state a generalization of 
Williams' decomposition of the excursion of Brownian motion at the maximum,
given in Chaumont \cite{chaumont94}. We need to
make a step out of the world of probability and consider $\sigma$-finite
measures instead of probability laws. Recall that $m_v=-\un{X}_v$ is the 
absolute value of the minimum before time $v$, and with our notations
$T_{m_v-}$ is the first time (and a.s.\ last before $v$) when $X$ attains
this value. Write
$$\begin{array}{cll}
\xg_s&=X_s & 0\leq s\leq T_{m_v-},\\
\xd_s&=m_v+ X_{s+T_{m_v-}}&  0\leq s\leq v- T_{m_v-}
\end{array} $$
for the pre- and post- minimum processes of $X$ before time $v$. Then by 
\cite{chaumont94},
\begin{lmm}\label{willdec}
One has the identity for
$\sigma$-finite measures
$$\int_0^{\infty}\d v P^v(\xg\in \d\omega, \xd\in \d \omega')=
\int_0^{\infty}\d x P^{(-x,\infty)}(\d \omega)\otimes \int_0^{\infty}\d u
N^{>u}(\d \omega'),$$
where $N^{>u}$ is the finite measure characterized by 
$N^{>u}(F(X))=N(F(X_s,0\leq s\leq u),\zeta(X)>u)$ for every non-negative
measurable $F$. 
This
in turn determines entirely the laws $P^v$ for $v>0$. 
\end{lmm}
Loosely speaking, if $v$ is ``random'' with ``law'' the Lebesgue measure 
on $(0,\infty)$, the pre- and
post- minimum processes are independent with respective ``laws''
$\int_0^{\infty}\d x P^{(-x,\infty)}(\d \omega)$ and $\int_0^{\infty}\d u
N^{>u}(\d\omega)$.  
As a consequence of this identity, we have that under $P^v$ for some
fixed $v>0$,
conditionally on $m_v$ and $T_{m_v-}=\tau$, the processes $\xg$ and $\xd$
are independent with respective laws $P^{\tau}_{0\downarrow -m_v}(\d \omega)$
and $(N^{> v-\tau}(1))^{-1}N^{> v-\tau}(\d \omega')$.

\begin{lmm}\label{condint}
Let $z>0$. Under the probability $P^1_{0\to -z}$, 
conditionally on $\tau_3-\tau_1=t$, the ranked sequence of 
lengths of the constancy intervals of the infimum process of 
$(X_{s+\tau_1},0\leq s\leq \tau_3-\tau_1)$ have the same law as 
$\Delta T_{[0,z]}$ given $T_z=t$ under $P$. 
\end{lmm}

\proof
We first condition by the value of $(m, \tau_3)$. Then 
by Lemma \ref{willdec} the path $\underleftarrow{X}$ has the law 
$P^{\tau_3}_{0\downarrow -m}$ 
of the first-passage bridge from $0$ to $-z$ with
lifetime $\tau_3$. Applying Lemma \ref{tatb} and the Markov property
we obtain that conditionally 
on $\tau_1$ the path $(X_{s+\tau_1}+m-z,0\leq s\leq \tau_3-\tau_1)$ is 
a first passage bridge ending at $-z$ at time $\tau_3-\tau_1$. Since it 
depends only on $\tau_3-\tau_1$, we have 
obtained the conditional distribution given $\tau_3-\tau_1$.
Hence, the sequence defined in the lemma's statement has the same conditional
law as the ranked lengths of 
the constancy intervals of the infimum process of such a first-passage bridge, 
that is, it has the same law as $\Delta T'_{[0,z]}$ given 
$T'_z=\tau_3-\tau_1$, with $T'$ as in the statement. \cq

The last lemma gives an explicit form for the law of the 
remaining length  $1-\tau_3+\tau_1$ under $P^1_{0\to -z}$.

\begin{lmm}\label{lastint}
One has 
$$P^1_{0\to -z}(1-\tau_3+\tau_1\in \d s)=\d s\frac{c_{\alpha}zq_z(1-s)}{
s^{1/\alpha}q_z(1)},$$
which is the law of a size-biased pick of the sequence $\Delta T_{[0,z]}$ 
given $T_z=1$ under $P$. 
\end{lmm}

\proof
By Lemma \ref{willdec}, 
if $s$ is ``distributed'' according to Lebesgue measure on $\R_+$, 
then under $P^s$, the processes $\underleftarrow{X}$ and 
$\underrightarrow{X}$ are independent with respective ``laws''
$\int_0^{\infty}\d x P^{(-x,\infty)}(\d\omega)$ and $\int_0^{\infty}\d u 
N^{>u}(\d\omega')$. 
Our first task is to disintegrate these laws 
to obtain a relation under $P^1_{0\to -z}$. Let $H$ and $H'$ be two 
continuous bounded functionals and $f$ be continuous with a compact support 
on $(0,\infty)$. Then, letting $T^{\omega}_{\cdot}=
\inf\{s\geq 0:\omega(s)<\cdot\}$,   
\begin{eqnarray*}
\lefteqn{\int_0^{\infty}\!\!\d s f(s)E^s[H(\underleftarrow{X})
H'(\underrightarrow{X})\, |\,
|X_s+z|<\eps]}\\
&=&\int_0^{\infty}\!\!\!\d x\int_0^{\infty}\!\!\!\d u \int\!\!\!
\int P^{(-x,\infty)}(\d\omega)
N^{>u}(\d\omega')f(T^{\omega}_x+u)H(\omega)H'(\omega')
\frac{\ind_{\{|z-x+\omega'(u)|<\eps\}}}{P(|X_{T^{\omega}_x+u}+z|<\eps)}\\
&=&\int_0^{\infty}\!\!\!\d u\int\!\! N^{>u}
(\d\omega')H'(\omega')\!
\int_0^{\infty}\!\!\!\d x \frac{\ind_{\{|z-x+\omega'(u)|<\eps\}}}{2\eps}\!\!
\int\! P^{(-x,\infty)}(\d\omega)\frac{f(T^{\omega}_x+u)H(\omega)}{(2\eps)^{-1}
P(|X_{T^{\omega}_x+u}+z|<\eps)}.
\end{eqnarray*}
The measure $(2\eps)^{-1}\ind_{\{|z-x+\omega'(u)|< \eps\}}\d x$ converges 
weakly as $\eps\to0$ to the Dirac mass at $z+\omega'_u$. Recall that the 
family of probability measures $P^{(-x,\infty)}$ is continuous as $x$
varies. Since $f$ has compact support, we can restrain $T^{\omega}_x+u$ to 
stay in a compact set. Then, the denominator in the last integral, 
which converges to $p_{T^{\omega}_x+u}(-z)$, remains
bounded and converges uniformly in $x$ and $u$. 
Then the boundedness of $H$ implies that the two last integrals
converge to 
$$\int P^{(-z-\omega'_u,\infty)}(\d\omega)\frac{f(T^{\omega}_{z+\omega'(u)}+u)
}{p_{T^{\omega}_x+u}(-z)}.$$
Now, the measure $N^{>u}$ is a finite measure, so the fact that
$u$ actually stays in a compact set and the fact that the two last 
integrals above remain bounded allow to apply the dominated convergence 
theorem to obtain 
\begin{eqnarray*}
\lefteqn{\int_0^{\infty}\d s f(s)P^s_{0\to-z}(H(\underleftarrow{X})
H'(\underrightarrow{X}))}\\
&=&\int_0^{\infty}\!\!\!\d u \int\!
N^{>u}(\d\omega')H'(\omega')\int\! P^{(-z-\omega'(u),\infty)}(\d\omega)
H(\omega)\frac{f(T_{z+\omega'(u)}(\omega)+u)}{p_{T^{\omega}_x+u}(-z)}
\end{eqnarray*}
Now we disintegrate this relation by taking $f(s)=(2\eps)^{-1}
\ind_{[1-\eps,1+\eps]}(s)$, so a similar argument as above gives that 
the left hand side converges to $P^1_{0\to-z}(H(\underleftarrow{X})
H'(\underrightarrow{X}))$ as 
$\eps\downarrow0$, whereas the right hand side is 
$$\int_0^{\infty}\!\!\! \d u \int\! N^{>u}(\d\omega')
H'(\omega')\int\! P^{(-z-\omega'(u),\infty)}(\d\omega)H(\omega)\frac{\ind_{
[1-\eps,1+\eps]}(T^{\omega}_{z+\omega'(u)}+u)}{2\eps p_{
T^{\omega}_{z+\omega'(u)}
+u}(-z)}.$$
The third integral may be rewritten as 
$$\frac{P(|T^{\omega}_{z+\omega'(u)}+u-1|< \eps)}{2\eps}
E^{(-z-\omega'(u),\infty)}\left[\left.\frac{H(\omega)}{
p_{T^{\omega}_{z+\omega'(u)}+u}(-z)}
\right||T^{\omega}_{z+\omega'(u)}+u-1|< \eps\right],$$
with a slightly improper writing (the $\omega$'s should not appear in the 
expectation, but we keep them to keep the distinction with the expectation 
with respect to $\omega'$).
Similar arguments as above imply that the limit we are looking for is 
$$P^1_{0\to-z}(H(\underleftarrow{X})H'(\underrightarrow{X}))=
p_1(-z)^{-1}\int_0^1\!\!\! \d u N^{>u}\left[H'(\omega') 
q_{z+\omega'(u)}(1-u)E^{1-u}_{0\downarrow -(z+\omega'(u))}[H(\omega)]\right].$$
This in turn completely determines the law of the bridge by a monotone class
argument. A careful application of the above identity thus gives
$$E^1_{0\to-z}[f(1-(\tau_3-\tau_1))]=p_1(-z)^{-1}\int_0^1 \d u
N^{>u}\left[q_{z+\omega'(u)}(1-u)
E^{1-u}_{0\downarrow-(z+\omega'(u))}[f(u+T^{\omega}_{\omega'(u)})]\right].$$
Applying Lemma \ref{tatb} to the rightmost expectation term, this is equal to 
\begin{eqnarray*}
\lefteqn{p_1(-z)^{-1}\int_0^1 \d u
N^{>u}\left[q_{z+\omega'(u)}(1-u)\int_0^{1-u}\d
v\frac{q_{\omega'(u)}(v)q_z(1-u-v)}{q_{\omega'(u)+z}(1-u)}f(u+v)\right]}\\
&=&p_1(-z)^{-1}\int_0^1 \d u\int_u^1\d s f(s) q_z(1-s)
N^{>u}\left[q_{\omega'(u)}(s-u)\right]\\
&=&zq_z(1)^{-1}\int_0^1 \d s f(s)q_z(1-s)\int_0^s \d u
N^{>u}\left[q_{\omega'(u)}(s-u)\right]
\end{eqnarray*}
It remains to compute the second integral. Using scaling identities for 
$N^{>u}$ and $q_x(s)$ we have
\begin{eqnarray*}
\int_0^s \d u
N^{>u}\left[q_{\omega'(u)}(s-u)\right]&=&
\int_0^1\d
rN^{>sr}\left[q_{\omega'(sr)}(s(1-r))\right]\\
&=&s^{-1/\alpha}\int_0^1 \d r
s^{1/\alpha}N^{>sr}\left[q_{s^{-1/\alpha}\omega'(sr)}(1-r)\right]\\
&=&s^{-1/\alpha}\int_0^1\d r
N^{>r}\left[q_{\omega'(r)}(1-r)\right].
\end{eqnarray*}
Finally, the integral in the right hand side
does not depend on $s$, we call it $c$ and obtain
$$E^1_{0\to-z}[f(1-(\tau_3-\tau_1))]=
\int_0^1\d sf(s)\frac{czq_z(1-s)}{
s^{1/\alpha}q_z(1)}.$$ 
So we
necessarily have
$c=c_{\alpha}$, and the claim follows. \cq

\noindent{\bf Proof of Proposition \ref{prpsemig}. }
The proof is now easily obtained by combining the last lemmas. 
Under $\bP^1_{0\to 0}$, conditionally on $Z^{(t)}_1=z$, the law of
the lengths of constancy intervals of $\un{V^{\tau_2}X}^{(t)}$ 
is obtained by adjoining the term $1-(\tau_3-\tau_1)$ to a sequence which, 
conditionally on $1-(\tau_3-\tau_1)=t$, has same law  
as $\Delta T_{[0,z]}$ given $T_z=1-t$ under $P$ (Lemma \ref{condint}). 
By Lemma \ref{lastint}, 
$1-(\tau_3-\tau_1)$ has itself the law of a size-biased pick from 
$\Delta T_{[0,z]}$ given $T_z=1$ under $P$, so Lemma \ref{permanform}
shows the whole sequence has the law of 
$\Delta T_{[0,z]}$ given $T_z=1$. Last, by Lemma \ref{decompX}, $Z^{(t)}_1$ 
has density $p_1^{(t)}(-z)\rho_1^{(t)}(z)p_1(0)^{-1}\d z$, entailing the
claim. \cq

\subsection{Proof of Theorem \ref{T2}}

To recover the dislocation measure of $F^{\natural}$, we use
the following 
variation of Lemma \ref{msprp} and 
\cite[Corollary 1]{MS03}. For details on size-biased versions of 
measures on $S$, 
see e.g.\ \cite{donjoy89}, which deals with probability measures, but the 
results we mention are easily extended to $\sigma$-finite measures. 

\begin{prp}\label{msdisl}
Let $(F(t),t\geq 0)$ be a ranked self-similar fragmentation with
characteristics $(\beta,0,\nu)$, $\beta\geq 0$. For every $t$, let $F_*(t)$ be 
a random size-biased permutation of the sequence $F(t)$ (defined on a possibly 
enlarged probability space). 
Let $G$ be a continuous bounded function
on the set of non-negative sequences with sum $\leq 1$, depending only on the 
first $I$ terms of the sequence, with support included in a set of the 
form $\{s_i\in[\eta,1-\eta], 1\leq i\leq I\}$. Then
$$\frac{1}{t}E[G(F_*(t)]\build\to_{t\downarrow 0}^{}\nu_*(G),$$ 
where $\nu_*$ is the 
size-biased version of $\nu$ characterized by 
$$\nu_*(G)=\int_S \nu(\d\bs)\sum_{j_1,\ldots,j_I}G(s_{j_1},\ldots,s_{j_I})
s_{j_1}\frac{s_{j_2}}{1-s_{j_1}}\ldots\frac{s_{j_I}}{1-s_{j_1}-\ldots-
s_{j_I}},$$
where the sum is on all possible distinct $j_1,\ldots,j_I$.
Moreover, $\nu$ can be recovered from $\nu_*$.
\end{prp}

\noindent{\bf Proof of Theorem \ref{T2}. }
Let $G$ be a function of the form $G(x)=f_1(x_1)\ldots f_I(x_k)$ for 
$x=(x_1,x_2,\ldots)$ and $\sum_i x_i\leq 1$, with 
$f_1,\ldots,f_I$ continuous bounded functions on $[0,1]$ that are null on 
a set of the form $[0,1]\setminus]\eta,1-\eta[$. Let 
$\Delta^* T_{[0,z]}$ be the sequence of the jumps of $T$ on the interval 
$[0,z]$, listed in size-biased order (which involves some 
enlargement of the probability space). Using Lemma \ref{permanform}, 
it is easy that 
$z\mapsto E[G(\Delta^* T_{[0, z]})|T_z=1]$ is a continuously differentiable 
function with derivative bounded by some $M>0$. Let also 
$F^{\natural}_*(t)$ be the sequence $F^{\natural}(t)$ listed in size-biased order. 
Now by Proposition \ref{prpsemig}, 
$$\bN^{(1)}\left[\frac{G(F^{\natural}_*(t))}{t}\right]=\frac{1}{t}\bE\left[
e^{-t^{\alpha}+tZ^{(t)}_1}
p_1(-Z^{(t)}_1)p_1(0)^{-1} \wt{E}\left[G(\Delta^* \wt{T}_{
[0,Z^{(t)}_1]})\left|
\wt{T}_{Z^{(t)}_1}=1\right.\right]\right],$$
where $\wt{T}$ is a copy of $T$ with law $\wt{E}$, independent of the marked 
process $X$. Consider a function $f(t,z)$ that is continuous in $t$ and $x$ 
and null at $(t,0)$ for every $t\geq 0$. Then the compensation formula applied the
subordinator $Z^{(t)}$ between times $0$ and $1$ gives
\begin{eqnarray}
\frac{1}{t}\bE[f(t,Z^{(t)}_1)]&=&\frac{1}{t}\int_0^1\d x
\int C_{\alpha}(1-e^{-ts})
s^{-\alpha-1}\d s
\bE[f(t,Z^{(t)}_x+s)-f(t,Z^{(t)}_x)]\nonumber\\
&\build\to_{t\to 0}^{}& C_{\alpha} \int_0^1\d x\int s^{-\alpha}\, \d s\, 
f(0,s)=C_{\alpha}\int s^{-\alpha} \d sf(0,s), \nonumber
\end{eqnarray}
as soon as we may justify the convergence above. 
Take 
$$f(t,z)=\exp(-t^{\alpha}+tz)p_1(-z)p_1(0)^{-1}
E[G(\Delta T^*_{[0,z]})|T_z=1],$$ 
then we have to check that 
$s^{-\alpha}\bE[|f(t,Z^{(t)}_x+s)-f(t,Z^{(t)}_x)|]$ is bounded independently 
on $x\in[0,1]$. By the 
hypotheses on $G$, it is again true that $z\mapsto f(t,z)$ is a 
continuously differentiable function with uniformly bounded derivative, 
when $t$ stays in a neighborhood of $0$. 
Hence the expectation above is bounded 
by $(M' s\wedge M'')s^{-\alpha}$ 
for some $M',M''>0$, which allows to apply the dominated 
convergence theorem. By Proposition \ref{msdisl}, 
we obtain, denoting by $\nu_{\natural}$ the dislocation measure of 
$F^{\natural}$, 
\begin{eqnarray*}
t^{-1} \bN^{(1)}[G(F^{\natural}_*(t))]&\build\to_{t\to0}^{}& \int_S 
\nu_{\natural}(\d\bs) 
G(s_{j_1},\ldots,
s_{j_I}) 
\sum_{j_1,\ldots,j_I}s_{j_1}\frac{s_{j_2}}{
1-s_{i_1}}\ldots \frac{s_{j_I}}{1-s_{j_1}-\ldots-s_{j_{I-1}}}\\
&=& C_\alpha \int_0^{\infty} \d s\frac{s^{-\alpha}p_1(-s)}{p_1(0)} 
E[G(\Delta T^*_{[0,s]})
|T_s=1],
\end{eqnarray*} 
allowing to conclude that $\nu_{\natural}=\nu_+$ 
with the same computations as in the proof of Lemma
\ref{rform}. 
\cq

\section{Asymptotics}\label{smasym}

In this section we discuss asymptotic results for $F^+$. 

\subsection{Small-time asymptotics}

\begin{prp}\label{smasymfp}
Let $Z$ be a non-negative stable ($\alpha-1$) random variable with Laplace
transform $E[\exp(-\lambda Z)]=\exp(-\alpha\lambda^{\alpha-1})$. 
Denote
by $\Delta_1,\Delta_2,\ldots$ the ranked jumps of $(T_x,0\leq x\leq Z)$,
where $T$ is as before the stable $1/\alpha$ subordinator, which is 
taken independent of $Z$. 
Then
$$t^{\alpha/(1-\alpha)}(F^+_2(t),F^+_3(t),\ldots)
\build\to_{t\to0+}^{d}(\Delta_1,\Delta_2,\ldots).$$
\end{prp}

We first need the

\begin{lmm}\label{cvgstm}
Let $Z^{(t)}_1$ have the law $\rho_1^{(t)}(s)\d s$ above, then 
$$t^{1/(1-\alpha)}Z^{(t)}_1\build\to_{t\to0+}^{d}Z,$$
where $Z$ is as above a stable variable with Laplace exponent
$\alpha\lambda^{\alpha-1}$. 
\end{lmm}

\proof
Recall that $Z^{(t)}$ is a subordinator with characteristic exponent given by 
$$\bE[e^{-\lambda Z^{(t)}_1}]=\exp\left(-
\int_0^{\infty}\frac{C_{\alpha}(1-e^{-tx})\d x}{
x^{\alpha+1}}(1-e^{-\lambda x})\right).$$ 
Therefore, evaluating the Laplace exponent at the point $t^{1/(1-\alpha)}
\lambda$, changing variables and using dominated convergence entails 
$$\bE[\exp(-\lambda t^{1/(1-\alpha)} Z^{(t)}_1)]\build\to_{t\to 0+}^{}
\exp\left(-\int_0^{\infty}\frac{C_{\alpha}
\d y}{y^{\alpha}}(1-e^{-\lambda y})\right).$$
Thus the convergence to some limiting $Z$. 
Using now the explicit value for $C_{\alpha}$, we see that the Laplace 
exponent of $Z$ has to be $\alpha\lambda^{\alpha-1}$, as claimed. \cq

The proof of Proposition \ref{smasymfp} follows the same lines as for 
Proposition 6 in \cite{mierfmoins}, 
so we will only sketch it. One first begins 
with proving that if ${\cal Z}$ 
is as in Lemma \ref{decompX} a random variable distributed according
to the law that has density $\rho^{(t)}_1(z)p^{(t)}_1(-z)\d z/p_1(0)$, then 
$t^{1/(1-\alpha)}{\cal Z}$ 
converges in law to $Z$. This is a consequence
of the preceding lemma, since as $t\to0$, $X^{(t)}$ converges to $X$, 
so one can write 
$$\bE[g(t^{1/(1-\alpha)}{\cal Z})]=
\bE[g(t^{1/(1-\alpha)}Z^{(t)}_1)p^{(t)}_1
(-Z^{(t)}_1)/p_1(0)],$$
where $Z^{(t)}_1$ is distributed as above. By Skorokhod's representation 
theorem, 
we may suppose that $t^{1/(1-\alpha)}Z^{(t)}_1$ converges a.s.\ to its limit
in law $Z$, So it remains to show that a.s.\ $p^{(t)}_1(-Z^{(t)}_1)\to
p_1(0)$  as $t\to0$ to apply dominated convergence, and this is done by 
recalling that $p_1^{(t)}(z)=e^{-t^{\alpha}-tz}p_1(z)$.
Then one reasons by induction
just as in \cite[Proposition 6]{mierfmoins}, using the explicit form of the 
semigroup of $F^+$. 

\subsection{Large-time asymptotics}

By a direct application of Theorem 3 in \cite{bertafrag02}, one gets
the large $t$ asymptotic behavior for $F^+$. Recall that the Gamma
law with parameter $a$ is the law with density proportional to 
$x^{a-1}e^{-x}$ on $\R_+$. The moments of this law are given, for 
$r>-a$, by 
$$\frac{1}{\Gamma(a)}
\int_0^{\infty} x^{r+a-1}e^{-x}\d x=\frac{\Gamma(a+r)}{\Gamma(a)}.$$
\begin{prp}\label{asymfp}
Define 
$$\rho_t(\d y)=\sum_{i=1}^{\infty}F_i(t)\delta_{t^{\alpha}F_i(t)}(\d y),$$
then $\rho_t$ is a probability measure that converges in law as $t\to\infty$ 
to the deterministic Gamma law with parameter 
$1-1/\alpha$.
\end{prp}

\proof
We know by \cite[Theorem 3]{bertafrag02} that $\rho_t$ converges to 
some probability $\rho_{\infty}$ that is characterized by its moments, 
$$\int_0^{\infty} y^{k/\alpha}\rho_{\infty}(\d y)=\frac{\alpha(k-1)!}{
\Phi'(0+)\Phi\left(\frac{1}{\alpha}\right)\ldots \Phi\left(\frac{k-1}{\alpha}
\right)}$$
for every $k\geq 1$, where $\Phi$ is the Laplace exponent of a subordinator
related to a tagged fragment of the process $F^+$. This exponent depends
only on the dislocation measure (and not the index), so it is the same 
as for $F_-$ in \cite{mierfmoins}. 
By taking the explicit value of $\Phi$ (Section 3.2 therein), 
we easily get
$$\int_0^{\infty}y^{k/\alpha}
\rho_{\infty}(\d y)=\left(\frac{\alpha \Gamma\left(
1+\frac{1}{\alpha}\right)}{\Gamma\left(\frac{1}{\alpha}\right)}\right)^k
\frac{\Gamma\left(1+\frac{k-1}{\alpha}\right)}{\Gamma\left(1-\frac{1}{\alpha}
\right)}=\frac{\Gamma\left(1+\frac{k-1}{\alpha}\right)}{\Gamma\left(1-\frac{1}{\alpha}
\right)}.$$
Replacing $k$ by $\alpha k$, one can recognize the moments of the 
Gamma law with the claimed parameter. \cq

\noindent{\bf Acknowledgments. }
Many thanks to Jean Bertoin for many precious comments on this work, 
and to Jean-Fran\c cois Le Gall for discussions related to the stable tree. 
Thanks also to an anonymous referee for a careful reading and 
very helpful comments that helped to consequently improve the presentation of
this work.

\end{document}